\documentclass{article}

\def\Z{{\rm Z}}
\def\N{{\rm N}}
\def\C{{\rm C}}

\def\R{{\rm R}}
\def\Q{{\rm Q}}

\begin{document}

Tsemo Aristide,

Coll\`ege Bor\'eal

Campus de Toronto

951, avenue Carlaw

Toronto ON M4K 3M2

tsemo58@yahoo.ca

\bigskip
\bigskip

\centerline{\bf Differentiable Categories, differentiable gerbes
and $G$-structures.}

\bigskip
\bigskip

\centerline{\bf Abstract.}

\bigskip

The theories of strings and $D$-branes have motivated the
development of non Abelian cohomology techniques in differential
geometry, on the purpose to find a geometric interpretation of
characteristic classes. The spaces studied here, like orbifolds
are not often smooth. In classical differential geometry, non
smooth spaces appear also naturally, for example in the theory of
foliations, the space of leaves can be an orbifold with
singularities. The scheme to study these structures is identical:
classical tools used in differential geometry, like connections
and curvatures are adapted. The purpose of this paper is to
present the notion of differentiable category which unifies all
these points of view. This enables us to provide a geometric
interpretation of $5$-characteristic classes, and to interpret
classical problems which appear in the theory of $G$-structures by
using gerbes.

\bigskip
\bigskip

\section{Introduction.}

\bigskip

Differential geometry is the study of the analytic properties of
topological spaces. Most of the main tools developed in this
theory are issued from calculus, and  a manifold does not have
singularities. Moduli spaces in differential geometry are rarely
smooth as show the space of orbits of the action of a compact Lie
group on a manifold, the space of leaves of a foliation, the
compactification of a space of curves,... Since these singular
structures arise naturally in differential geometry, it is normal
to try to study them by using methods created in the smooth case.
An example of such a method is the theory of orbifolds created by
Satake [37], which enables to study the action of finite groups on
manifolds  which may have fixed points (see Audin [4]), foliations
with bundle-like metrics (see Molino and Pierrot [33]), strings
theory (see Chen and Ruan [11]), quotient of compact affine
manifolds, in particular quotient of flat affine spaces forms are
studied by many authors (see Long and Reid [20], Ratcliffe and
Tschantz [36]).... We can also quote other theories such as the
theory of foliages which studies the differential geometry of the
space of leaves of a foliation (see Molino [31]), the study of
homogeneous $(X,H)$-manifolds see (Goldman [15]).

\medskip

One of the main goal of this paper is to propose the theory of
differentiable categories to unify the generalizations of
classical differential geometry mentioned above: a differentiable
category is a category whose objects are differentiable manifolds,
and the morphisms between its objects are differentiable maps.
This point of view enables to handle also new situations like
generalized orbifolds such as the orbit space of the action of a
compact Lie group, or the space of leaves of a foliation endowed
with a bundle-like metric. (see Molino and Pierrot [33]).

\medskip

In mathematics, the classification problem is the cornerstone on
which relies every theory $T$; this is performed by assigning to
objects which occur in $T$ simpler invariants which enable to
describe them completely: for example, the genus of a closed
surface. The scheme usually followed in classical differential
geometry is to define objects and invariants locally, and glue
them with partitions of unity. Local invariants in the theory of
differentiable categories are more difficult to study, since even
when there exists a topology, neighborhoods of different points
are not always isomorphic, for example the notion of a frames
bundle is not straightforward defined since the dimension of the
objects in a differentiable category may vary. This situation is
analog to algebraic geometry, and we intensively use the machinery
developed by Grothendieck and his students in this context (see
Giraud [13], [14]). In fact sheaves of categories and gerbes are
nowadays intensively studied by differential geometers (see
Brylinski [8], Brylinski and McLaughlin [9], [10]): the functional
action in classical mechanic which describes the motion of a point
is expressed by using a connection on a principal bundle. In the
purpose to unify all existing fundamental strengths, physicists
have defined strings and branes theories. The functional action
which describes the motion of a string is defined by a gerbe, and
one expects that a good notion of $n$-gerbes will enable to handle
branes theories. In fact sheaves of categories in this context are
examples of differentiable categories. A currently very active
research topic is the adaptation of tools defined in classical
geometry like connections on principal bundles to these objects.

\medskip

We start this paper by studying the differential geometry of
differentiable categories without using Grothendieck topologies.
We define the notions of principal bundles, which are torsors
whose fibers are Lie groups,  the tangent space  of a
differentiable category and its DeRham cohomology. In this setting
we introduce connections forms and distributions and study their
holonomy.

In modern geometry, global objects are constructed by gluing local
objects. For example, a manifold is obtained by gluing open
subsets of a vector space, schemes in algebraic geometry are
defined by gluing spectrums of commutative rings. Algebraic
geometers have remarked that in many situations, the transition
functions that are used to glue objects do not verify the Chasles
relation. This has motivated descent theory which is presented in
the setting of categories theory by Giraud [13]. We study
differential descent; or equivalently descent in the theory of
differentiable categories. This is an adaptation of the analysis
situs of Giraud; we introduce differentiable fibered principal
functors and their connective structures. Recall that the notion
of connective structure has been introduced by Brylinski [8] in
the context of gerbes on manifolds to provide a geometric
interpretation of characteristic classes.

\medskip

The local analysis intensively used in differential geometry
relies on the existence of neighborhoods of points. This is
achieved in this context by differentiable Grothendieck
topologies: examples of Grothendieck topologies are defined on
orbifolds, generalized orbifolds, foliages,... the CechDeRham
complex is then used to study cohomology. Chen and Ruan [11] have
defined a new cohomology theory for orbifolds to understand
mathematical strings theory. We adapt to generalized orbifolds
this new cohomology theory.

With the notion of Grothendieck topologies defined, we can study
sheaves of categories and gerbes in the theory of differentiable
categories. The first example of such a construction can be
obtained by gerbes defined on  the Grothendieck topology
associated to an orbifold. Lupercio and Uribe [21] have provided
such a construction by using groupoids. One of the fundamental
example of a differentiable gerbe is the canonical gerbe defined
on a compact simple Lie group $H$ (see Brylinski [8]). The
classifying cocycle of this gerbe is the canonical $3$-cohomology
class defined by the Killing form. Medina and Revoy [26], [27]
have classified Lie groups endowed with a non degenerated
bi-invariant scalar product which also defines a canonical
$3$-form. Remark that these Lie groups are not always compact and
are even contractible when they are nilpotent and simply
connected. The theory of lattices in Lie groups presented by
Raghunathan [35] and the Leray-Serre spectral sequence enable us
to construct fundamental examples of gerbes on compact manifolds
which are the quotient of a nilpotent Lie group by a lattice.

The notion of Grothendieck topology of a differentiable category
enables us to construct the curving, and the curvature of a
connective structure on a differentiable principal gerbe. We also
define the holonomy form which is used to study functional action
on loop spaces.

\medskip

An approach of the study of the differential geometry of a gerbe
can be done by using right invariant distributions defined in the
thesis of Molino [29]. We outline how to a gerbe defined on a
manifold one can associate an invariant distribution which enables
to construct the holonomy around curves.

\medskip

In the last part of the paper, we study sequences of fibered
categories. An example of such a construction has been done by
Brylinski and McLaughlin [9] to provide a geometric representation
of the Pontryagin  class of degree $4$. To a principal gerbe, we
associate a $2$-sequence of fibered categories which must be an
example  of a $U(1)$ $3$-gerbe (recall that the notion of
$3$-gerbe is not well-understood yet). We associate to such a
$2$-sequence of fibered categories  a $5$-integral cohomology
class.

\medskip

This new tools for differential geometers can be used to tackle
well-known problems in differential geometry. A $G$-structure is a
reduction of the bundle of jets defined on a manifold. This theory
has been intensively studied in the seventies (see Molino [32] and
the thesis of Albert [1],  Medina [25], Nguiffo-Boyom [34]). We
can associate to a manifold a sheaf of categories which represents
the geometric obstruction to the existence of a $G$-structure.

\bigskip

{\bf Plan.}

\bigskip

1. Introduction.

\medskip

2. Notation.

\medskip

3. Basic definitions and examples.

3.1 Orbifolds $(X,H)$-manifolds and differentiable categories.

3.2 Actions of compact Lie groups and differentiable categories.

3.3 Foliages and differentiable categories.

3.4 Projective presented manifolds.

\medskip

4. Differentiable fibered categories.

4.1 Connections on differentiable bundles torsors.

4.2 Differentiable tensors of a differentiable category.

4.3 Frames bundle and differentiable categories.

4.4 Differentiable descent and connections in fibered categories.

\medskip

5. Grothendieck topologies in differentiable categories.

5.1 Grothendieck topologies and cohomology of differentiable
categories.

5.2 Chen-Ruan cohomology for generalized orbifolds.

\medskip

6. Sheaf of categories and gerbes in differentiable categories.

6.1 The classifying cocycle of a gerbe.

\bigskip

7. Examples of sheaves of categories and gerbes.

7.1 Gerbes and $G$-structures.

7.2 Gerbes and invariant scalar products on Lie groups.

7.3 A sheaf of categories on an orbifold with singularities.

\medskip

8. Differential geometry of sheaves of categories.

8.1 Induced bundles.

8.2 Reduction to the motivating example.

8.3 Holonomy and functor on loop spaces.

8.4 Canonical relations associated to the connective structure on
a gerbe.

8.5 Uniform distribution and gerbes.

\medskip

9. Sequences of fibered categories in differentiable categories.

9.1 $4$-cocycles and sequences of fibered categories

References

\section{ Notations.}

Let $C$ be a category which has a final object $e$, and $I$ a
small set relatively to a given universe (see  [3] SGA 4 p. 4). In
fact the cardinality used throughout this paper are numerable.
Since we are studying differentiable manifolds, we want our spaces
to be at least paracompact, to insure existence of partitions of
unity, one of the main tools used in differential geometry to show
the existence of global objects.

Consider a small family $(X_i)_{i\in I}$ of objects of $C$. We
denote by $X_{i_1...i_n}$ the fiber product (if it exists) of the
finite subset $\{X_{i_1},...,X_{i_n}\}$ of $(X_i)_{i\in I}$ over
$e$.

Let $P$ be a presheaf of categories defined over $C$. For every
objects $e_i, {e'}_i\in P(X_{i_1})$, and  a map $u:e_i\rightarrow
{e'}_i$, we denote by $e^{i_2...i_n}_i$ and by $u^{i_2...i_n}$ the
respective  restrictions of $e_i$ and $u$ to $U_{i_1...i_n}$.

\section{Basic Definitions and examples.}

The differentiable manifolds used in this paper are $C^{\infty}$,
and finite dimensional.

\medskip

{\bf Definition 3.1.}

A differentiable category $C$ is a category such that:

- every element $X$ of the class of objects of $C$  is a
differentiable manifold,

- every morphism of $C$  is a differentiable map.

\medskip

{\bf Examples.}

\medskip

The category $Diff$, whose objects are finite dimensional
differentiable manifolds, and such that the set of morphisms
$Hom_{Diff}(M,N)$ between two differentiable manifolds $M$ and $N$
is the set of differentiable maps between $M$ and $N$, is a
differentiable category. Remark that this category is not small
relatively to an universe $U$ which contains the set of real
numbers, but is $U$-small (see [3] SGA4 p. 5).

\medskip

Let $N$ be a manifold, and $C_N$ the category whose objects are
open subsets of $N$. The  morphisms between objects of $C_N$ are
the canonical imbeddings. The category $C_N$ is endowed with the
structure of a differentiable category, for which each open subset
of $N$ is endowed with the differentiable structure inherited from
$N$.

\subsection{Orbifolds, $(X,H)$-manifolds and differentiable
categories.}

The theory of orbifolds has been introduced by Satake (see Satake
[37], Chen and Ruan[11]).  Orbifolds appear in different branches
of mathematics, like strings theory, foliations theory: the
singular foliation defined by the adherence of the leaves of a
foliation endowed with a bundle-like metric can define an orbifold
(see Molino and Pierrot [33] p. 208).

\medskip

{\bf Definition 3.1.1.}

 An $n$-dimensional orbifold $N$ (see also
Chen and Ruan [11] definition 2.1), is a separated topological
space $N$, such that:

 - for every element $x\in N$, there exists an open
subset $U_x$ of $N$,

- an open subset $V_x$  of  ${\R}^n$, a finite group of
diffeomorphisms $\Gamma_x$ of $V_x$, an element $\hat x\in V_x$

such that for every element $\gamma_x$ in $\Gamma_x, \gamma_x(\hat
x)=\hat x$.

-There exists a continuous map $\phi_x:V_x\rightarrow U_x$, such
that $\phi_x(\hat x)=x$, for each $y\in V_x$, and $\gamma_x$ in
$\Gamma_x$, $\phi_x(\gamma_x(y))=\phi_x(y)$, and the induced
morphism $V_x/\Gamma_x\rightarrow U_x$ is an homeomorphism.

The triple $(V_x,\phi_x,\Gamma_x)$ is called an orbifold chart.

 We suppose that the following condition is satisfied:

 Let $(V_x,\phi_x,\Gamma_x)$ and $(V_y,\phi_y,\Gamma_y)$ be two orbifolds charts.
 We denote by $p_x:V_x\times_NV_y\rightarrow V_x$ the canonical projection.

We suppose that there exists an equivariant diffeomorphism in
respect of $\Gamma_x$ and $\Gamma_y$:

$$
\phi_{xy}:p_y(V_x\times_NV_y)\rightarrow p_x(V_x\times_NV_y)
$$

such that:

$$
{\phi_y}_{\mid p_y(V_x\times_NV_y)}={\phi_x}_{\mid
p_x(V_x\times_NV_y)}\circ\phi_{xy}.
$$

The fact that the morphism $\phi_{xy}$ is equivariant is
equivalent to saying that for every element $\gamma_y$ in
$\Gamma_y$, there exists an element $\Phi_{xy}(\gamma_y)$ in
$\Gamma_x$ such that:

 $$
 \phi_{xy}\circ
\gamma_y=\Phi_{xy}(\gamma_y)\circ\phi_{xy}.
$$.

Remark that the fiber product $V_x\times_NV_y$ is not necessarily
a manifold. This can be illustrated by the following example:
consider the quotient $N$, of the real line ${\R}$, by the map
$x\rightarrow -x$, the fiber product ${\R}\times_N{\R}$ is the
union of two non parallel lines in ${\R}^2$; but
$p_x(V_x\times_NV_y)$ is an open subset of $V_x$.

The maps $\phi_x\circ\phi_{xy}\circ{\phi_{yz}}_{\mid
p_z(V_x\times_NV_y\times_NV_z)}$ and $\phi_x\circ
{\phi_{xz}}_{p_z(\mid V_x\times_NV_y\times_NV_z)}$ are equal.
Since $\Gamma_x$ is finite, there exists an element $c_{xyz}$ in
$\Gamma_x$ such that:

$$
\phi_{xy}\circ{\phi_{yz}}_{\mid
p_z(V_x\times_NV_y\times_NV_z)}=c_{xyz}{\phi_{xz}}_{p_z(V_x\times_NV_y\times_NV_z)}.
$$

Let $\gamma_z$ be an element of $\Gamma_z$, we have:

$$
\phi_{xy}\circ\phi_{yz}\circ{\gamma_z}_{\mid
p_z(V_x\times_NV_y\times_NV_z)}=\Phi_{xy}(\Phi_{yz}(\gamma_z))\circ
c_{xyz}{\phi_{xz}}_{p_z(\mid V_x\times_NV_y\times_NV_z)}.
$$

We also have:

$$
c_{xyz}{\phi_{xz}}\circ{\gamma_z}_{p_z(V_x\times_NV_y\times_NV_z)}=
c_{xyz}\circ\Phi_{xz}(\gamma_z)\circ
{\phi_{xz}}_{p_z(V_x\times_NV_y\times_NV_z)}
$$

This implies that:

$$
\Phi_{xy}\circ\Phi_{yz}=c_{xyz}\circ\Phi_{xz}\circ{c_{xyz}}^{-1}.
$$

\medskip

An example of an orbifold is the quotient of a manifold by a
finite subgroup of its group of diffeomorphisms. We are going to
characterize some orbifolds arising from the theory of affine
manifolds. (See also  Long and Reid [20];  Ratcliffe, and Tschantz
[36]).

Let $H$ be a Lie group which acts differentiably and transitively
on the manifold $X$, we first recall some basic facts of the
theory of $(X,H)$ manifolds (See Goldman [15]).

A $(X,H)$-manifold, is a differentiable manifold, endowed with an
open covering $(U_i)_{i\in I}$ such that:

 - For each $i\in I$, there
exists a differentiable map $d_i:U_i\rightarrow X$, such that
$d_i:U_i\rightarrow d_i(U_i)$ is a diffeomorphism.

- The transition function ${d_j\circ {d_i}^{-1}}_{\mid d_i(U_i\cap
U_j)}:d_i(U_i\cap U_j)\rightarrow d_j(U_i\cap U_j)$ coincides with
the restriction of the action of an element of $H$ on $d_i(U_i\cap
U_j)$.

An $(X,H)$-map $f:N\rightarrow N'$ between the $(X,H)$-manifolds
$N$ and $N'$, is a differentiable map which preserves their
$(X,H)$ structures.

Examples of $(X,H)$ structures are $n$-dimensional affine
manifolds: here $H$ is the group $Aff({\R}^n)$ of affine
transformations of ${\R}^n$, and $X={\R}^n$; $n$-dimensional
projective manifolds where $X$ is the real projective space
$P^{n}({\R})$, and $H$ is $PGl(n,{\R})$ the group of projective
transformations. A $({\R}^n,Aff({\R}^n))$-automorphism is called
an affine transformation.

We can show the following result relating orbifolds to affine
manifolds, and to projective manifolds:

\medskip

{\bf Proposition 3.1.1.}

{\it Let $N$ be an affine manifold, and $\Gamma$ a finite group of
affine transformations of $N$, such that the set $C$ of elements
of $N$, such that for every element $u$ of $C$, there exists a non
trivial element $\gamma$ of $\Gamma$ such that $\gamma(u)=u$, is
finite; moreover we suppose that every element of $C$ is fixed by
every element of $\Gamma$. Let $p:N\rightarrow N/\Gamma$ be the
canonical projection. The blowing-up (not in the classical sense)
of $N/\Gamma$ at $p(C)$ is a projective manifold.}

\medskip

{\bf Proof.}
 First we are going to blow-up the action of
$\Gamma$.

Let $u$ be an element of $C$, and $U$ an affine chart around $u$.
Thus $U-\{u\}$ can be identified with a ball without the origin.
Consider the submanifold ${P'}^n$ of ${\R}^n\times P^{n-1}{\R}$
defined by the equations:

 $$
 (x_1,...,x_n,[X_1,...,X_n])\in{P'}^n \Longleftrightarrow
x_iX_j-x_jX_i=0.
$$

 There exists a projection
$p':{P'}^n\rightarrow {\R}^n$, the restriction of the canonical
${\R}^n\times P^{n-1}{\R}\rightarrow {\R}^n$. The blowing up of
$N$ at $u$ is the operation which replaces $U$ by ${p'}^{-1}(U)$
(see McDuff-Salamon [23] p. 233-235, See also Tsemo [39]). We can
cover $P^{n-1}{\R}$ by two open affine subsets $U_1$ and
$U_2$,which are the trivializations of the   the ${\R}$-line
bundle ${P'}^n$ over $P^{n-1}{\R}$. The coordinates change of
these trivializations is the map:

$$
u_{12}:U_1\cap U_2\times {\R}\longrightarrow U_1\cap U_2\times
{\R}
$$

$$
(x,y)\longrightarrow (x,-y)
$$

 Thus the imbedding maps:

 $$
 u_i:i=1,2:U_i\times{\R}\rightarrow P^{n+1}{\R}
 $$

 $$
 (x,y)\longrightarrow [x,y]
 $$

 defines a projective structure  around $u$ which can be glued
 with the affine atlas of $N-C$ to obtain a projective structure
 on the blowing-up $\hat N$, of $N$.

 We can identify the restriction of the action of the element $\gamma$ of $\Gamma$
on $U$, to a linear map $A_{\gamma}$, and extends it to a map
$A'_{\gamma}$ of ${P'}^n$ defined by
$A'_{\gamma}(x_1,...,x_n[X_1,...,X_n])=(A_{\gamma}(x_1,...,x_n),A_{\gamma}([X_1,...,X_n])$.
We thus obtain a free action of $\Gamma$ on $\hat N$ by projective
maps. The quotient of $\hat N$ by this action is a projective
manifold $N"$, obtained from $N/\Gamma$, by replacing a
neighborhood of every element $p(u)$, $u\in C$ by the quotient of
${p'}^{-1}(U)$ by $\Gamma$. We also say that $N"$ is a blowing-up
(not in the classical sense) of $N/\Gamma$ $\bullet$

\medskip

We associate to a orbifold $N$ the differentiable category ${C_N}$
defined as follows:

An object of $C_N$ is a triple $(M,\phi_M,\Gamma_M)$ where $M$ is
a manifold, $\Gamma_M$ a finite group of diffeomorphisms of $M$,
and $\phi_M:M\rightarrow N$, a continuous map such that for every
element $\gamma_M$ in $\Gamma_M$, for every element $x\in M$,
$\phi_M(\gamma_M(x))=\phi_M(x)$, and the induced map
$M/\Gamma_M\rightarrow N$ is a local homeomorphism. For every
$y=\phi_M(x)$ in $N$, there exists a chart of the orbifold
$(V_x,\phi_x,\Gamma_x)$ around $y=\phi_x(x)$ (see definition
3.1.1) such that:  if $p^M:V_x\times_NM\rightarrow M$, and
$p_x:V_x\times_NM\rightarrow V_x$ are the natural projections,
 there exists an equivariant local diffeomorphism
$\phi^M_{V_x}:p^M(V_x\times_NM)\rightarrow p_x(V_x\times_NM)$ such
that $\phi_x\circ \phi^M_{V_x}={\phi_M}_{\mid p^M(V_x\times_NM)}$.
In particular a chart of $M$ is an object of ${C_N}$.

A morphism between the objects $(M,\phi_M,\Gamma_M)$, and
$(M',\phi_{M'},\Gamma_{M'})$, is an equivariant differentiable map
$\phi:(M,\Gamma_M)\rightarrow (M',\Gamma_{M'})$ such that
$\phi_M=\phi_{M'}\circ \phi$.

\subsection{ Actions of compact Lie groups and differentiable
categories.}

 Let $M$ be a finite dimensional manifold, and $G$ a   compact Lie group which
 acts effectively on $M$. This is equivalent to saying that if an element of $G$ fixes every
element of $M$, it is the identity. We denote by $N$ the quotient
of $M$ by $G$. The orbits of $G$ are submanifolds, and when the
action is free, a well-known elementary result implies that $N$ is
a manifold (see Audin [4] p. 13-19). In the general situation $N$
is an orbifold with singularities. This can be seen by using the
slice theorem of Koszul that we recall now:

\medskip

{\bf Theorem 3.2.1 Koszul [18].}

{\it Let $G$ be a compact Lie group which acts effectively and
differentiably on the manifold $M$. Let $u$ be an element of $M$.
Denote by $G_u$ the subgroup of $G$ which fixes  $u$, there exists
an invariant neighborhood $U$ of $u$ which is isomorphic to a
neighborhood of the zero section in the quotient of $G\times V$ by
$G_u$, where $V$ is the quotient of the tangent space $T_uM$ by
its subspace tangent to the orbit.}

 \medskip

Thus, the slice theorem allows to construct an open covering of
$M/G$ whose elements are quotient of  open subsets of a vector
space by the action of a compact Lie group (take a transversal to
the zero section in theorem 3.2.1).

\medskip

Let $N$ be the quotient space of $M$ by $G$. We associate to the
action of $G$ on $M$ the following differentiable category $C_N$
defined as follows:

 An object of ${C_N}$ is a triple $(P,H,\phi_P)$
where $P$ is a manifold endowed with an effective  action of a
compact Lie group $H$, such that there exists  a local equivariant
diffeomorphism $\phi_P:(P,H)\rightarrow (M,G)$ such that the
induced map $P/H\rightarrow N=M/G$ is a local homeomorphism.

 A morphism $f$ between the
objects $(P,H,\phi_P)$ and $(P',H',\phi_{P'})$ is defined by  an
equivariant differentiable map $f:P\rightarrow P'$ such that
$\phi_{P'}\circ f=\phi_{P}$.

\medskip

The previous construction can be generalized in the following
setting:

\medskip

{\bf Definition 3.2.1.}

 Let $N$ be a separated topological space, A generalized orbifold
 on $N$ is defined by the following data:

 For every element $u\in N$, there exists a manifold $M_u$, a
 compact Lie group $H_u$ which acts differentiably on $M_u$, and a continuous map
 $\phi_u:M_u\rightarrow N$ whose image contains $u$; such that for
 every $h_u$ in $H_u$, for every $x$ in $M_u$,
 $\phi_u(h_u(x))=\phi_u(x)$, and the induced map $M_u/H_u\rightarrow
 N$ is a local homeomorphism. The triple $(M_u,H_u,\phi_u)$ is called a
 chart of the generalized orbifold.

 Let $(M_u,H_u,\phi_u)$ and $(M_{v},H_{v},\phi_{v})$ be two
 charts. Denote by $p_u:M_u\times_NM_{v}\rightarrow M_u$ the
 canonical projection. There exists a local equivariant
 diffeomorphism $\phi_{uv}:p_v(M_u\times_NM_v)\rightarrow
 p_u(M_u\times_NM_v)$ such that ${\phi_v}_{\mid
 p_v(M_u\times_NM_u)}=\phi_u\circ \phi_{uv}$.

\bigskip

Moduli spaces appear in different domains of differential
geometry; many of them can be endowed with the structure of a
differentiable category $C$. For example, consider the
differentiable category whose class of objects is the class of
isomorphism classes of hyperbolic surfaces of genus $h$, $h$
fixed. Let $[X_h]$ be the class of the surface of genus $h$,
$X_h$. The differentiable structure of $[X_h]$ is the
differentiable structure of one element picked in the class of
$[X_h]$, for example $X_h$ itself. The set of morphisms
$Hom([X_h],[X_{h'}])$ is the set of hyperbolic maps between
$X^1_h$ and $X^2_{h'}$, where $X^1_h$ and $X^2_{h'}$ are the
respective representants   picked in the classes of $[X_h]$ and
$[X_{h'}]$ to define the structure of the differentiable category.

\subsection{ Foliages and differentiable categories.}

Let $N$ be a $n$-dimensional manifold, a foliation ${\cal F}$ on
$N$ of codimension $q$ is defined by an atlas
$(U_i,\phi_i:U_i\rightarrow {\R}^p\times {\R}^q )_{i\in I}$, such
that $\phi_i\circ {\phi_j}^{-1}_{\mid \phi_j(U_i\cap
U_j)}(x,y)=(u_{ij}(x,y),v_{ij}(y))$. This is equivalent to define
a partition of $N$ by immersed manifolds of dimension $p$ called
the leaves. In this situation we say that the couple $(N,{\cal
F})$ is a foliated manifold. One of the main important problem in
foliation theory is the study of the topology of the space of
leaves, which is not always endowed with the structure of a
manifold. For example consider the quotient $T^2$ of ${\R}^2$ by
the group $\Gamma$ generated two translations $t_{e_1}$ and
$t_{e_2}$ whose directions $e_1$ and $e_2$ are independent
vectors. Let $\theta$ be an irrational integer; the foliation of
${\R}^2$ by affine lines parallel to $e_1+\theta e_2$ defines on
$T^2$ a foliation for which every leaf is dense. Thus the space of
leaves of this foliation is not separated.

\medskip

In [31] Molino has introduced the notion of foliage to study these
situations which can be interpreted with differentiable
categories:

\medskip

{\bf Definitions 3.3.1.}

Two foliated manifolds $(N_1,{\cal F}_1)$ and $(N_2,{\cal F}_2)$
are transversally equivalent if and only if there exists a
foliated manifold $(\hat N,\hat{\cal F})$, two submersions
$\pi_i,i=1,2: \hat N\rightarrow N_i$ such that the leaves of
$\hat{\cal F}$ are the preimages of the leaves of ${\cal F}_i$ by
$\pi_i$.

Let $N$ be a topological space a foliage on $N$ is a
differentiable category $C_N$ whose objects are quadruples
$(U,V,{\cal F},\pi)$, where $U$ is an open subset of $N$,
$(V,{\cal F})$ a foliated manifold. We assume that the space of
leaves of ${\cal F}$ is $U$ and $\pi:V\rightarrow U$ is the
natural projection.

Let $(U',V',{\cal F}',\pi')$ another object of $C_N$, we denote by
$p_V:V\times_NV'\rightarrow V$ the natural projection. We assume
the quadruples $(\pi(p_V(V\times_NV')),p_V(V\times_NV'),{\cal
F}_{\mid p_V(V\times_NV')},\pi_{\mid p_V(V\times_NV')})$ and
$(\pi(p_{V'}(V\times_NV')),p_{V'}(V\times_NV'),{\cal F}'_{\mid
p_{V'}(V\times_NV')},\pi'_{\mid p_{V'}(V\times_NV')})$ are objects
of $C_N$ and transversally equivalent, where ${\cal F}_{\mid
p_V(V\times_NV')}$ is the restriction of ${\cal F}$ to
$p_V(V\times_NV')$.

Finally we suppose that for every element $y\in N$, there exists
an object $(U,V,{\cal F},\pi)$ of $C_N$, such that $y\in U$.

A morphism between the objects $(V,U,{\cal F},\pi)$ and
$(V',U',{\cal F}',\pi')$ of $C_N$ is a differentiable map
$\phi:V\rightarrow V'$ such that $\pi=\pi'\circ \phi$.

\subsection{ Projective presented manifolds.}

In order to study a generalized equivalence Cartan problem, Molino
(see Molino [32]) has studied projective presented manifolds which
are examples of differentiable categories:

\medskip

{\bf Definition 3.4.1.}

A projective presented manifold is a small differentiable category
whose class of objects is a projective system of manifolds
$(V_i,\pi^i_j:V_j\rightarrow V_i)_{i\in I}$. The following
conditions need to be satisfied:

The maps $\pi^i_j$ are submersions,

Let $\hat V$ be the topological  projective limit of the family
$(V_i,\pi^i_j)_{i\in I}$, and $(x_i)_{i\in I}, x_i\in V_i$ an
element of $\hat V$. Let $P^i:\hat V\rightarrow V_i$ which
associates to $(x_i)_{i\in I}$ the element $x_i$ in $V_i$. For
each $i$, there exists an open neighborhood $U_i$ of $x_i$ in
$V_i$, a map $c_i:U_i\rightarrow \hat V$, such that $P^i\circ
c_i=Id_{U_i}$, and $P^j\circ c_i$ is differentiable.

Let $(x_i)_{i\in I}$ be an element of $\hat V$, there exists
$i_0$, and a neighborhood $U_{i_0}$ of $x_{i_0}$ in $V_{i_0}$ such
that for every $i>i_0$, for every $y\in U_{i_0}$,
${\pi^{i_0}_{i}}^{-1}(y)$ is connected in $U_i$.

\medskip

{\bf Definition 3.4.2.}

 Let $C$ and $C'$ be two differentiable
categories; a differentiable morphism between $C$ and $C'$ is
defined by a functor $F:C\rightarrow C'$, such that for every
object $X$ of $C$, there exists a differentiable map
$h^F_X:X\rightarrow F(X)$, such that for every morphism
$f:X\rightarrow X'$ in $C$, the following square is commutative:

$$
 \matrix{X &{\buildrel {f}\over{\longrightarrow}} & X'\cr \downarrow
h^F_X &&\downarrow h^F_{X'}\cr F(X) & {\buildrel
{F(f)}\over{\longrightarrow}} & F(X')}
$$

\medskip

We can suppose that the categories $C$ and $C'$ are imbedded in
$Diff$ the category  of differentiable manifolds. In this setting,
a differentiable functor is a morphism between the identity
functor of $C$, and a functor $F:C\rightarrow Diff$ whose image is
contained in $C'$, moreover   for every object $X$, the map
$h^F_X:X\rightarrow F(X)$ which defines the morphism of functors
is a differentiable map.

\medskip

{\bf Examples.}

Let $f:M\rightarrow N$ be a differentiable map, $f$ can be viewed
as a differentiable functor $F:C\rightarrow C'$  where the unique
object of $C$ is $M$ and the unique object of $C'$ is $N$. We
suppose that the only morphisms in $C$ and $C'$
 are the identities. The functor $F$ assigns $N$ to $M$, and
 $h^F_M=f$.

 \medskip

 Suppose that a Lie group $G$ acts differentially on $M$, and $N$,
  we define $C$, to be the differentiable category which has $M$ has a unique
 object, and such that $Hom_C(M,M)$ is the image of the map
 $G\rightarrow Diff(M)$ which defines the action. Similarly, we define $C'$ to be the category whose
 unique object is $N$  and such that $Hom_{C'}(N,N)$ is the image of the map
 $G\rightarrow Diff(N)$.
 Let $\phi$ be an endomorphism of $G$, each $\phi$-equivariant map
 $f:M\rightarrow N$; that is a map such that for each $g\in G$,
 $f\circ g=\phi(g)\circ f$ defines a differentiable functor $F$ between $C$
 and $C'$, such that $F$ assigns $N$ to $M$,  $h^F_M=f$, and $F(g)=\phi(g)$.

\section{ Differentiable fibered categories.}

To study the differentiable structure of differentiable
categories, we are going to use the theory of fibered categories.
On this purpose, we recall the following facts adapted to our
setting:

\medskip

{\bf Definition 4.1.}

Let $F:P\rightarrow C$ be a differentiable functor, and
$f:x\rightarrow y$ a map of $C$. Let $x',z'$ be two objects  of
the fiber of $x$, and $y'$ an object of the fiber of $y$. Denote
by $Hom_{f}(z',y')$  the subset of the set of morphisms
$Hom_P(z',y')$ such that for every element $l\in Hom_{f}(z',y')$,
$F(l)=f$.

 A morphism
$f':x'\rightarrow y'$ is Cartesian, if and only if the map
$Hom_{Id_{x}}(z',x')\rightarrow Hom_{f}(z',y')$ which assigns to
$h$ the map $f'\circ h$ is bijective for every $z'$ in the fiber
of $x$.

\medskip

{\bf Definition 4.2.}

 A differentiable bundle functor  $F:P\rightarrow C$ is a Cartesian
 functor which satisfies the following conditions:

 - The fiber of an object $x$ of $C$ has a unique element $p_x$.

 - For every object $x$ of $C$, there exists a Lie group $H_x$ such that
   the canonical projection
$p_x\rightarrow x$ defines on $p_x$ the structure of a total space
of a $H_x$-principal bundle, whose
 base space is $x$.  Morphisms between objects of $P$ are
 morphisms between  principal differentiable bundles.

 \medskip

 If the group $H_x$ is independent of $x$, we say that $F:P\rightarrow C$ is a $H$-principal
 differentiable bundle functor.

 \medskip

 {\bf Example.}

 Let $H$ be a compact group which acts on the manifold $N$ by
 diffeomorphisms. We have attached a differentiable category $C_N$
 to this action (see p. 8). Let $(P,H_P,\phi_P)$, an object of $C_N$, we can
 construct the principal $H_P$-bundle $P_{H_P}$ which is the quotient of
 $P\times H_P$ by the diagonal action of $H_P$. Let
 $f:(P,H)\rightarrow (P',H')$ a morphism in $C_N$ which is induced by a morphism
 $l_{H,H'}:H\rightarrow H'$ such that for every elements $h$ in $H$,
 and $p$ in $P$, $f(hp)=l_{H,H'}(h)f(p)$ since by definition $f$ is an equivariant map. We
 deduce a morphism $\psi_{H,H'}(f):P\times H\rightarrow P'\times H'$ which sends
 $(p,h)$ to $(f(p),l_{H,H'}(h))$. For every $h_0\in H$, we have:

 $$
 \psi_{H,H'}(f)(h_0p,h_0h)=(f(h_0p),l_{H,H'}(h_0h))=l_{H,H'}(h_0)\psi(p,h).
 $$

Thus the morphism $\psi_{H,H'}(f)$ induces a morphism
$\psi'_{H,H'}(f):P_{H}\rightarrow P'_{H'_{}}$.

 We deduce the existence of a differentiable category $PC_N$
 whose class of objects are the bundles $P_{H_P}$, and a
 differentiable bundle functor $F_N:PC_N\rightarrow C_N$ which
 sends the object $P_{H_P}$ to $P$. The Cartesian map above
 $f$ is $\psi'_{H,H'}(f)$.

\medskip

{\bf Definition 4.3.}

Let $F:P\rightarrow C$ be a differentiable bundle functor, and
$H:C'\rightarrow C$, a morphism between differentiable categories.
The pull-back of $F$ by $H$ is the differentiable bundle functor
$F':P'\rightarrow C'$ defined as follows:

Let $X'$ be an object of $C'$, $h^H_{X'}:X'\rightarrow F(X')$ the
map which defines $H$; let $p_{F(X')}$ the object of the fiber of
$F(X')$ for $F$, and $h_{p(F(X'))}:p_{F(X')}\rightarrow F(X')$ the
bundle map. The fiber of $X'$ is the fiber product of the maps
$h^H_{X'}$ and $h_{p(F(X'))}$.

\subsection{ Connection on differentiable bundle functors.}

In this part we are going to study connections on differentiable
bundles functors. First we recall the notion of connection on a
principal bundle (See Lichnerowicz [19] p. 56, McDuff and Salamon
[23] p. 207-209).

 Let $H$ be a Lie group whose Lie algebra is denoted by ${\cal H}$,
 and $p:P\rightarrow N$ a $H$-principal bundle over the $n$-dimensional
 manifold $N$, for every element $A\in {\cal H}$,
 we denote $A^*$ the vector field defined on $P$ by the formula:

 $$
 A^*(x)=lim_{t\rightarrow 0}{d\over{dt}}xexp(tA), x\in P.
 $$

 A connection  defined on the $H$-principal bundle
 $p:P\rightarrow N$, is a $1$-form $\theta:P\rightarrow {\cal H}$
 which verifies the following conditions:

 Let $A^*$ be the fundamental vector field defined by $A\in {\cal H}$,  $i_{A^*}\theta=A$.

 For every element $h\in H$,  $h^*\theta=Ad(h^{-1})\theta$.

A connection is also defined by a distribution on $P$ transverse
to the fibers and invariant by $H$, whose rank is $n$ the
dimension of $N$. To a connection form $\theta$, the distribution
associated is: $\Theta_x=\{u\in TP_x,\theta(x)=0\}$.

 The curvature of $\theta $ is the ${\cal H}$-valued $2$-form
 on $P$ defined by:
  $\Omega=d\theta+{1\over 2}[\theta,\theta]$.

We adapt now this definition to differentiable categories:

  \medskip

  {\bf Definition 4.1.1.}

  A connection  on the principal bundle functor
  $p:P\rightarrow C$ is defined  by a connection form $\theta_X$ on the
  principal $H_X$-bundle $p_X\rightarrow X$ of the fiber of $X$, such that for a map $h:p_X\rightarrow p_Y$
  (necessarily Cartesian), the distribution defined by the kernel of
   $h^*(\theta_Y)$ is the distribution which defines
  the connection form of $\theta_X$.

  \medskip

  {\bf Example.}

  Consider the interval $I=]-1,1[$ of ${\R}$, and $N$ the orbifold which is the quotient
    of $I$ by the symmetry $h:x\rightarrow -x$, we associate to this orbifold the differentiable
    category $C_N$ whose class of objects contains only $I$, and the set of
    morphisms of $I$, $Hom_{C_N}(I,I)=\{Id_I,h\}$, remark that this is not the canonical
    differentiable category associated to an orbifold defined at p. 8. The real  $1$-form
    $\alpha=xdx$ is invariant by $h$, thus defines a connection on
    the trivial bundle functor $P\rightarrow C_N$ in circles over $C_N$ as follows:
    let $C^1$ be the circle, $P$ is the category which unique object is $e_I=I\times C^1$.
     The unique non trivial morphism of $e_I$ is the map $h'$
     defined by $h'(x,y)=(-x,y)$. Let $(u,v)$ be a vector tangent to
     $(x,y)\in I\times C^1$ we set ${\theta_I}_{(x,y)}(u,v)=\alpha_x(u)+v=xu+v$.

\bigskip

{ \bf Definition: Holonomy of a connection of a principal  bundle
functor 4.1.2.}

\medskip

 Consider $C_{I}$, the canonical differentiable category defined
 on the interval by its structure of manifold, and $F:P\rightarrow C$ a $H$-principal bundle functor,
 endowed with a connection form $\theta$ and
 $L:C_{I}\rightarrow C$ a differentiable functor.
 The pull-back of $F$ and $\theta$ by $L$ is a principal bundle over the
 interval endowed with a connection form whose holonomy map is the holonomy of $F:P\rightarrow C$, around
 $L$.

\subsection{ Differentiable tensors of a  differentiable
category.}

In this section, we are going to associate to a differentiable
category $C$, principal bundles functors which allow to define
tensor fields. Such a theory is obviously known for manifolds. It
has also been developed in the category of orbifolds see (Chen and
Ruan [11]), and for foliages (see Molino [31]).

\medskip

{\bf Definition 4.2.1.}

Let $C$ be a differentiable category, the differentiable tangent
bundle of $C$ is the differentiable category $T(C)$ defined as
follows: the elements of the class of objects of $T(C)$ are
tangent spaces $T(X)$, where  $X$ is an object of $C$. A map
between $T(X)$ and $T(Y)$, is a map $T(h):T(X)\rightarrow T(Y)$
induced by a morphism $h:X\rightarrow Y$ in $C$.

A differentiable functor $F:C\rightarrow C'$, induces a tangent
functor $T(F):T(C)\rightarrow T(C')$ defined as follows: let $X$
be an object of $C$, the map $h^F_X:X\rightarrow F(X)$ induces the
tangent map $T(h^F_X):T(X)\rightarrow T(F(X))$ which defines the
tangent functor.

\medskip

The differentiable category of $p$-forms of $C$, $\Lambda^p(C)$,
is the category whose class of objects is the class whose elements
are $\Lambda^p T(X)$, where $X$ is an object in $C$. A morphism
between the objects $\Lambda^pT(X)$ and $\Lambda^pT(Y)$ is a map
of the form $\Lambda^p(h)$ where $h:X\rightarrow Y$ is a morphism
in $C$.

A differentiable $p$-form is a functor
$\alpha:\Lambda^p(C)\rightarrow C_{\R}$,  where $C_{\R}$ is
endowed with the structure of a differentiable category which as a
unique object: the real line ${\R}$, and such that the
endomorphisms of ${\R}$ are differentiable maps of ${\R}$. We
deduce from the definition of a differentiable functor that the
following condition is satisfied:  let $\alpha_X$ be a $p$-form,
for every map $f:X\rightarrow Y$, there exists a diffeomorphism
$\alpha(f)$ of ${\R}$ such that the following square is
commutative:

$$
\matrix{\Lambda^pX &{\buildrel {d^pf}\over{\longrightarrow}}&
\Lambda^p Y\cr \downarrow \alpha_X && \downarrow \alpha_Y\cr {\R}
& {\buildrel{\alpha(f)}\over{\longrightarrow}} & {\R}}
$$

Let $f:X\rightarrow Y$ be a morphism in $C$, we don't assume that
$\alpha(f)$ is the identity of ${\R}$. Thus $\alpha_X$ is not
necessarily the pull-back of $\alpha_Y$ by $f$. We denote by
$\Lambda^p_{Id}(C)$ the set of $p$-forms such that for every map
$f$ in $C$, $\alpha(f)$ is the identity.

\medskip

{\bf Definition-Proposition 4.2.2.}

{\it Let $C$ be a differentiable category, and $\alpha$ a $p$-form
defined on $C$, there exists a functor $d:\Lambda^p(C)\rightarrow
\Lambda^{p+1}(C)$, the differential such that $d\circ d=0$.}

\medskip

{\bf Proof.}
 Let $\alpha$ be a $p$-form defined on $C$, for each
object $X$ of $C$, $\alpha_X$ is a $p$-form, we can define
$d\alpha_X$ the differential of $\alpha_X$. Let $f:X\rightarrow Y$
be a morphism in $\Lambda^p(C)$, we define
$d(\alpha)(f)=d(\alpha(f))$$\bullet$

\medskip

{\bf Examples.}

\medskip

Let $N$ be a manifold, and $H$ a Lie group which acts
differentially on $N$. Consider the differential category $C_N^H$
whose unique object is $N$, and such that the set of endomorphisms
of $N$ in $C^H_N$ is the image of the map $H\rightarrow Diff(N)$
which defines the action. Let $\chi$ be a character of $H$, we can
define $\Lambda^p_{H,\chi}(N)$ to be the set of $p$-forms on
$C^H_N$ such that for every element $\alpha\in
\Lambda^p_{H,\chi}(N)$, the following square is commutative:

$$
\matrix{\Lambda^pN &{\buildrel {d^pf}\over{\longrightarrow}}&
\Lambda^p N\cr \downarrow \alpha_X && \downarrow \alpha_Y\cr {\R}
& {\buildrel{\chi(f)}\over{\longrightarrow}} & {\R}}
$$

In particular if $\chi$ is the trivial character, we obtain the
space of $H$-invariant $p$-forms, and the equivariant cohomology.

\medskip

Let $C_N$ be the differentiable category associated to a foliage.
The set of $\Lambda^p_{Id}C_N$ forms on $C_N$ is the set of basic
forms (see Molino [31]).

\subsection{ Frames bundle and differentiable categories.}

Let $C$ be a differentiable category, we cannot always define the
bundle of linear frames, since two objects of $C$ do not have
necessarily the same dimension. Suppose that every objects $C$ has
dimension $n$. We can define the set of vector frames $V(C)$ as
follows: Let $X$ be an object of $C$, and $x$ be an element of
$X$, we denote by $V(C)_x$ the set of linear maps
$u:{\R}^n\rightarrow T_xX$, where $T_xX$ is the tangent space of
$x$. We can thus define the vector bundle $V(C)(X)$ over $X$ whose
fiber at $x$ is $V(C)(X)_x$. Let $f:X\rightarrow Y$ be a
differentiable map, and $u\in V(C)(X)_x$, the linear map
$df_x\circ u$ is a vector frame of $V(C)(Y)_{f(x)}$. We have thus
define the category of vector frames of $C$. Since the morphisms
in $C$ are not necessarily local diffeomorphisms, we cannot assume
that the elements of $V(C)_x$ are isomorphisms.

\medskip

Let $C^H_N$ be the differentiable category associated to the
action of a compact Lie group $H$ on  $N$ (see p. 9). Let
$(P,H_P,\phi_P)$ be an object of $C_N$. We can define the vector
space $TCN_x$, the quotient of the tangent space at $x$,  $TP_x$
of $P$, by the image of the infinitesimal action at $x$ of $H_P$.
This space is called the tangent space at $x$ of the action.
Remark that the dimension $dim(TCN_x)$ of $TCN_x$ depends only of
$\phi_P(x)$. But this dimension can vary if $x$ varies in $P$.

 We can define the differentiable category of linear frames $L(C_N)$ of $C$. For each
object $X$ of $C_N$, there exists fibration $L(C_N)(X)\rightarrow
X$ such that for every element  $x\in X$, $L(C_N)(X)_x$ is the set
of linear isomorphisms ${\R}^{dim(TCN_x)}\rightarrow TCN_x$.

\medskip

{\bf Proposition 4.3.1.}

{\it Suppose that  the dimension of $TCN(P)_x$ does not depend of
$P$, then there exists a connection on the frames bundle $L(C_N)$,
of the differentiable category $C_N$ defined by the differentiable
category defined above.}

\medskip

{\bf Proof.}
 Let $\theta$ be a connection on the frames bundle
$L(C_N)(N)$ of $TCN(N)$ invariant by the action of $H$. Let
$(P,H_P,\phi_P)$ be an object of $C_N$. Since the dimension of
$TCN(P)_x$ does not depend neither of $P$, nor of $x$ in $P$, the
pull-back of $\theta$ by $\phi_P$ defines a connection form on
$L(C_N)(P)$$\bullet$

\medskip

{\bf Definition 4.3.1.}

Let $C_N$ be the differentiable category $C_N$ defined by the
action of the compact Lie group $H$ on $N$. We suppose that the
dimension of $TCN_x$ does not depend of $x$. Let $(X,H_X,\phi_X)$
be an element of $C_N$, and $\alpha_X$ the fundamental $1$-form of
the bundle $L(C_N)(X)$. It is the form ${\R}^{dim(TCN)}$-valued
form defined by:

$$
{\alpha_X}_u(v)=u^{-1}(dp_X(v))
$$

 where $u$ is an element of $L(C_N)(X)_x$, $v$ element of the
 tangent space of $L(C_N)(X)$ at $u$, and
 $p_X:L(C_N)(X)\rightarrow X$ the canonical projection. The family
 of $1$-forms $(\alpha_X)_{X\in C_N}$ defines an invariant form on
 $C_N$.

\subsection{Differentiable descent and connection in fibered
categories.}

\bigskip

In this part we are going to study the notions of connection and
holonomy on differentiable fibered categories.

Let us recall some facts on the analysis situs in differentiable
categories (see Giraud [13]):

\bigskip

Let $F:P\rightarrow C$ be a Cartesian functor, a clivage is a
family $L$ of morphisms  of $P$ such that:

every element in $L$ is cartesian,

for every morphism $f:x\rightarrow y$ in $C$, and $y'\in P_{y}$,
there exists a unique morphism $f'\in L$ whose target is $y'$ and
such that $F(f')=f$. A clivage is a scindage if and only if it is
stable by composition of maps. A clivage is the analog of a
reduction in differential geometry.

\medskip

Let $l:x\rightarrow y$ be a map in $C$, and $L$ a clivage. The
clivage $L$ and $l$ induce a functor $l^*:P_{y}\rightarrow P_x$
defined as follows:  The image of the object $z\in P_{y}$ is the
source of the unique Cartesian map $c_{l}(z):{l}^*(z)\rightarrow
z$ in $L$ over $l$.

Consider two maps  $l$ and $m$, such that the target of $l$ is the
source of $m$, there exists a natural transformation:

$$
c_{l,m}:{(m\circ l)}^*\rightarrow {l}^*\circ {m}^*
$$

which satisfies the relation:

$$
c_{ml}\circ c_{l,m}=c_{l}\circ c_{m}
$$

(See also Giraud [13] p.3).

 Let $p:P\rightarrow C$ be a Cartesian
functor between differentiable categories. We assume that there
exists a Lie group $H$ such that for every object $X$ of $C$,
every object $e_X$ in the fiber of $X$ is endowed with the
structure of an $H$-space; i.e the group $H$ acts on the right and
freely on $e_X$.

There exists a projection $p:e_X\rightarrow X$, such that for
every $h\in H$, and $x\in e_X$, $p(xh)=p(x)$.

A morphism $f:e_X\rightarrow e_{X'}$ in $C$ is a differentiable
map $f$ such that for every element $h\in H$, we have $f\circ
h=h\circ f$.

Let $f$ be an endomorphism of $e_X$, and $x\in e_X$. We denote by
$u(x)$ the element of $H$ such that $f(x)=xu(x)$. For every
element $h\in H$, we have $f(xh)=(xh)u(xh)=f(x)h=xu(x)h$. This
implies that:

$$
u(xh)=h^{-1}u(x)h.
$$

We suppose that there exists a principal $H$-bundle functor
$Aut(P)\rightarrow C$, such that for every object $e_X$ in the
fiber of $X\in C$, there exists a canonical isomorphism
$Aut(P)(X)\rightarrow End(e_X)$, which is natural in respect of
morphisms between objects.

Let $A$ be an element of ${\cal H}$ the Lie algebra of $H$, for
every object $e_X$, we can define the vector field:

$$
{d\over{dt}}_{t=0}xexp(tA).
$$
which is a fundamental vector field. This allows to identify
${\cal H}$ with a subbundle of the tangent space $T{e_X}_x$ of
$e_X$.

\medskip

Let us start by a motivating example. Let:

$$
1\longrightarrow H\rightarrow L'\rightarrow L\rightarrow 1
$$

be an exact sequence of Lie groups. Consider a principal
$L$-bundle $p:P\rightarrow N$ over the manifold $N$. The
obstruction to extend the structural group $L$, to $L'$, is
defined by a sheaf of categories $C_H$ defined as follows: for
every open subset $U$ of $N$, $C_H(U)$ is the category whose
objects are $L'$-principal bundles over $U$ whose quotient by $H$
is the restriction of $p$ to $U$. Morphisms between objects of
$C_H(U)$ are morphisms of $L'$-bundles which induce the identity
on the restriction of $p$ to $U$.

Let ${\cal L}$ and ${\cal L}'$ be the respective Lie algebras of
$L$ and $L'$. We know the definition of a connection form $\theta$
on $p$, and we want to generalize this definition. A natural way
is to take for each object $e_U\in C_H(U)$ a connection
$\alpha_U$, such that the composition of $\alpha_U$ with the
natural projection ${\cal L'}\rightarrow {\cal L}$ descends to the
restriction $\theta_U$ of $\theta$ to $U$. The choice of
$\alpha_U$ is not canonical since it is not necessarily preserved
by every automorphism $h$ of $e_U$.

Remark that  the form:

$$
\alpha_h=h^*(\alpha_U)-\alpha_U=Ad(h^{-1})(\alpha_U)-\alpha_U+h^{-1}dh
$$

is a ${\cal H}$-valued form.

\medskip

This motivates the following definition (compare with  Brylinski
[8] p. 206, and with Breen and Messing [7]):

\medskip

{\bf Definition 4.4.1.}

Let  $p:P\rightarrow C$ be a $H$-principal fibered category, and
${\cal H}$ the Lie algebra of $H$. A connective structure on $C$
is a map which assigns to every object $e_U$ of $P_U$, an affine
space $Co(e_U)$ such that:

 \medskip

The vector space of $Co(e_U)$ is the set of ${\cal H}$-forms
$\Omega^1(U,{\cal H})$.

\medskip

  For every  morphisms
$h':U'\rightarrow U"$, $h:U\rightarrow U'$, and for every object
$e_U$ in the fiber of $U$, $e_{U'}$ in the fiber of $U'$ and
$e_{U"}$ in the fiber of $U"$, there exists a morphism:

$$
h_*:Co(e_{U})\rightarrow Co(e_{U'})
$$

which is compatible with composition: $(h'h)_*={h'}_*h_*$.

There exists a morphism:

$$
u_h:h^*(Co(e_{U'}))\longrightarrow Co(h^*(e_{U'}))
$$

such that the following square is commutative:

$$
\matrix{{h}^*({h'}^*Co(e_{U"}))&
{\buildrel{u_{h'}}\over{\rightarrow}}& {h}^*Co({h'}^*e_{U"}) &
{\buildrel{u_{h}}\over{\rightarrow}}& Co({h}^*{h'}^*e_{U"})\cr
\downarrow \alpha_{h',h}^* &&&& \downarrow {{c_{h,h'}}^{-1}}_*\cr
(h'h)^*Co(e_{U"})
&&{\buildrel{u_{h'h}}\over{\longrightarrow}}&&Co((h'h)^*e_{U"})}
$$

where the morphisms  $c_{h,h'}$ is the morphism defined by a
morphism  in the analysis situs (see p. 16), and $\alpha_{h',h}$
the canonical isomorphism of torsors.

\bigskip

Let $u:e_U\rightarrow  e'_{U'}$ be a Cartesian morphism above
$h:U\rightarrow U'$, we have the compatibility diagram:

$$
\matrix{h^*Co(e_{U}) &{\buildrel{u_*}\over{\longrightarrow}}&
h^*Co(e'_{U'})\cr\downarrow u_h &&\downarrow u_h\cr Co(h^*e_U)
&{\buildrel{u_*}\over{\longrightarrow}}& Co(h^*(e'_{U'})).}
$$

\bigskip

There exists an action of $Aut_U(e_U)$ on $Co(e_U)$ such that for
every element $h$ of $Aut(e_U)$, and every element $\theta$ in
$Co(e_U)$. We have the relation:

$$
h_*(\theta)=h.(\theta)+h^{-1}dh.
$$

 And for every element $\alpha\in \Omega^1(U,{\cal H})$, we have:

 $$
 h_*(\theta+\alpha)=h_*(\theta)+Ad(h^{-1})(\alpha).
 $$

\medskip

An alternative definition of connective structure can be done by
considering torsors $Co(e_U)$ whose vector space is the space of
closed ${\cal H}$-valued $1$-forms if ${\cal H}$ is commutative.

\bigskip

{\bf A fundamental relation.}

\medskip

 Suppose now that $F:P\rightarrow C$ is a differentiable fibered
 category, consider a clivage $L$. For each objects $X$
 of $C$, and $X'$ in the fiber of $X$, consider
 a morphism $u_{XY}:Y\rightarrow X$, and its lift to
 a Cartesian morphism $u_{X'Y'}:Y'\rightarrow X'$ of $L$.

Let $\alpha$ be a connective structure defined on this
differentiable fibered bundle, we denote by $\alpha_{Y'}$ an
element of $Co(Y')$ , and by $\alpha_{X'Y'}$ the $1$-form such
that $\alpha_{X'}=\alpha_{X'Y'}+{u_{X'Y'}}_*(\alpha_{Y'})$ we
have:

$$
{u_{X'Y'}}_*(\alpha_{Y'Z'})-\alpha_{X'Z'}+\alpha_{X'Y'}=
$$

$$
={u_{X'Y'}}_*(\alpha_{Y'}-{u_{Y'Z'}}_*(\alpha_{Z'}))-(\alpha_{X'}-{u_{X'Z'}}_*(\alpha_{Z'}))
+(\alpha_{X'}-{u_{X'Y'}}_*(\alpha_{Y'}))
$$

$$
={u_{X'Z'}}_*(\alpha_{Z'})-{u_{X'Y'}}_*{u_{Y'Z'}}_*(\alpha_{Z'})
$$

Since $F:P\rightarrow C$ is a fibered category, there exists a
morphism $c_{X',Y',Z'}$ such that
$u_{X'Y'}u_{Y'Z'}=u_{X'Z'}c_{X',Y',Z'}$, we deduce that:

$$
{u_{X'Y'}}_*(\alpha_{Y'Z'})-\alpha_{X'Z'}+(\alpha_{X'Y'})=
$$

$$
{u_{X'Z'}}_*(\alpha_{Z'}-{c_{X',Y',Z'}}_*(\alpha_{Z'}))
$$

$$
={u_{X',Z'}}_*((\alpha_{Z'})-{c_{X',Y',Z'}}.(\alpha_{Z'})-{c_{X',Y',Z'}}^{-1}dc_{X',Y',Z'})
$$

\section{ Grothendieck topologies in differentiable categories.}

We have studied differentiable categories without emphasizing on
the global topology. This can be achieved by using the notion of
differentiable Grothendieck topology (see [3] S.G.A 4-1; p. 219;
or Giraud [14]).

\medskip

{\bf Definitions 5.1.}

Let $C$ be a differentiable category, a sieve $R$ in $C$ is a
subclass $R$ of the class of objects of $C$ such that if $U$ is an
object of $R$, and $V\rightarrow U$ is a morphism in $C$, then $V$
is in $R$.

A Grothendieck topology on $C$ is defined by assigning to each
object $U$ of $C$ a non empty family of sieves $J(U)$ of the
category over $U$, $C\uparrow U$ such that the following
conditions are satisfied:

For every morphism $h:V\rightarrow U$, and every sieve $R\in
J(U)$, the pull-back sieve $R^h$ is in $J(V)$.

A sieve $R$ of $C\uparrow U$ is in $J(U)$ if and only if for every
map $h:V\rightarrow U$, $R^h\in J(V)$.

\medskip

{\bf Examples.}

\medskip

 An example of a Grothendieck topology can be defined as follows:
 Let $N$ be a topological space, and $C_N$ the category whose
 objects are open subsets, and whose maps are canonical imbeddings between
 open subsets. For an
 open subset $U$, an element of $J(U)$ is a family of open subsets
 $(U_i)_{i\in I}$ of $U$ such that $\bigcup_{i\in I}U_i=U$. This
 topology is often called the small site.

\medskip

Let $N$ be a generalized orbifold (see definition 3.2.1). We can
define on $N$ the following Grothendieck topology:

A covering of an open subset $U$ of $N$, is a family of objects
$(P_i,H_i,\phi_i)_{i\in I}$ which is $U$-jointly surjective. This
equivalent to saying that $\bigcup_{i\in I}\phi_i(P_i)=U$. In
particular if for every object $(P,H_P,\phi_P)$ in $C_N$,  the
groups $H_P$ is discrete, we obtain a Grothendieck topology on
orbifolds.

\medskip

Consider the space of hyperbolic surfaces of a given genus $h$.
Each of this surface can be cut in pants. The hyperbolic length of
the boundaries cycles of these pants are the Fenschel-Nielsen
coordinates which identify the set of isomorphic classes of
hyperbolic surfaces of genus $h$ to a cell. (See [6] X. Buff and
al p.13-15).

\medskip

{\bf Definition 5.2.}

Let $(C,J)$ be a category endowed with a Grothendieck topology, we
suppose that $C$ has a final object $e$. A global covering of $C$
is a cover of $e$, that is an element of $J(e)$.

\medskip

{\bf Definition 5.3.}

A presheaf defined on the differentiable category $C$, is a
contravariant functor $F$, from $C$ to the category of sets.

A sheaf is a presheaf which satisfies $1$-descent in respect to
any sieve $R$ in $J(U)$. This is equivalent to saying that for
every object $U$ of $C$, and every sieve $R$ in $J(U)$, the
natural map:

$$
F(U)\rightarrow lim_{V\rightarrow U\in R}F(V)
$$
 is bijective.

\subsection{ Grothendieck topologies and cohomology of
differentiable categories.}

The cohomology of orbifolds is studied in algebraic geometry and
symplectic geometry, since orbifolds arise as phase spaces in
theoretical physics. Grothendieck and his collaborators (see
S.G.A. 4 II, p.16) have defined Cech cohomology in Grothendieck
sites. We shall apply this point of view to generalized orbifolds.
We shall also generalize  Chen and Ruan cohomology of orbifolds
(see [11]) to generalized orbifolds.

\medskip

 Let $J_N$ be the Grothendieck topology associated to the
generalized orbifold $N$. We can define the presheaf $\Omega^p_N$,
such that for each object $e=(P,H_P,\phi_P)$ of $C_N$,
$\Omega^p_N$ is the vector space of $p$-differentiable forms
invariant by $H_P$ defined on $P$ (see also p. 14). If
$h:e\rightarrow e'$ is a morphism in $C_N$, the restriction is
defined by the pull-back of differentiable forms.

Consider a covering $(U_i,H_i,\phi_i)_{i\in I}$ of $N$. We cannot
defined the classical Cech resolution, since the differentiable
category $C_N$ associated to $N$ is not necessarily stable fiber
products.
 Let
 $(U_{i_1},H_{i_1},\phi_{i_1}),...,(U_{i_n},H_{i_n},\phi_{i_n})$, be
 objects of $C_N$,
 $\phi_{i_1}(U_{i_1...i_n})$ is a manifold, We can
defined the bi-graded complex
$\Omega^l_N(\phi_{i_1}(U_{i_1...i_p}))$ endowed with two
derivations: the Cech-derivation and the canonical derivation of
differentiable forms. We denote by
$H^{*,*}_{(U_i,H_i,\phi_i)_{i\in I}}(N)$ the induced bigraded
cohomology groups.

We say that the covering $(U'_{i'},H'_{i'},\phi_{i'})_{i'\in I'}$
is finer than the covering $(U_i,H_i,\phi_i)_{i\in I}$, if and
only if for every $i'\in I'$, there exists $i\in I$  such that
$\phi_{i'}(U'_{i'})\subset \phi_i(U_i)$. This relation defines an
inductive system on the set of coverings, the inductive limit of
$H^{*,*}_{(U_i,H_i,\phi_i)_{i\in I}}(N)$ is the Cech-DeRham
cohomology of the generalized orbifold.

\subsection{ Chen-Ruan cohomology for generalized orbifolds.}

Suppose that the generalized orbifold $N$ is compact. We are going
to adapt the cohomology theory defined by Chen and Ruan [11] for
orbifolds. Firstly we recall the following construction in Chen
and Ruan (page 6-7): let $N$ be an orbifold, $(U_x,H_x,\phi_x)$ a
local chart at $x$, define $\hat N$ to be the set whose elements
are $(x,(h_x))$, where $(h_x)$ is the conjugacy class of the
element $h_x$ of $H_x$. Remark that $\hat N$ is well-defined
despite the use of local charts. The orbifold $\hat N$ is not
necessarily connected. Its connected components are called twisted
sectors (Chen and Ruan p.8). There exists a natural surjection
$p:\hat N\rightarrow N$, the connected components of elements of
$p^{-1}(U_x)$ can be parameterized by the set of conjugacy classes
$(h_x)$, $h_x$ in $H_x$.  Suppose that the orbifold is endowed
with a pseudo-complex structure, which defines a representation
$\rho_{H_x}:H_x\rightarrow Gl(n,{\C})$ ($n=dim_{{\C}}N$). For
every element $h_x$ in $H_x$, $\rho_{H_x}(h_x)$ depends only of
the conjugacy class  $(h_x)$ of $h_x$ in $H_x$, they define
$i_{x,(h_x)}=-{i\over{2\pi}}Log(det(\rho_{H_x}(h_x)))$. This
enables Chen and Ruan to define the orbifold $d$-cohomology group:

$$
H^d(X)=\oplus H^{d-2i_{(h)}}(X_{(h)}).
$$

\medskip

Let $N$ be a generalized compact orbifold, we can find a finite
cover $(U_i,H_i,\phi_i)$ for the Grothendieck topology, such that
each open subset $U_i$ is defined by the slice theorem (see
theorem 3.2.1), this is equivalent to saying that $U_i$ is the
quotient $H_i\times_{H'_i}V_i$ by $H_i$ where $H'_i$ is the
stabilizer of an element $x_i$ of $U_i$, $V_i$ is the quotient of
the tangent space $TU_i$ at $x_i$, by the image of the
infinitesimal action of $H_i$ at $x_i$ (see Audin [4] p. 15). Let
$C_i$ be an open subset of $V_i$ invariant by $H_i$. Then
$(C_i,H_i,\phi'_i)$ is a chart of the generalized orbifold, where
$\phi'_i:C_i\rightarrow C_i/H_i$ is the canonical projection. Thus
for every element $x$ in $N$, there exists a chart
$(U_x,H_x,\phi_x)$, $x'$ in $U_x$ such that $\phi_x(x')=x$, and
$H_x(x')=x'$. We are going to consider only this type of charts in
the sequel. The existence of such charts is related to the
definition of holonomy of singular foliations. See Molino and
Pierrot [33] p. 208, for the definition of the holonomy a
foliation defined by the action of a compact Lie group, or Debord
[12].

Let $H$ be a closed subgroup of $H_x$, we denote by $(H)$ the
conjugacy class of $H$ in $H_x$. Let $y$ be an element of $U_x$,
and $H^y_x$ the subgroup of $H_x$ which fixes $y$. We say that $y$
and $y'$ have the same type if and only if $(H^y_x)=(H^{y'}_x)$.
Let $(U_y,H_y,\phi_y)$ be a chart such that $H_y(y)=y$, denote by
$\lambda_y:H_y\rightarrow H_x$ the morphism induced by the
transition function $\phi_{xy}$. We   suppose that the stabilizer
of $\phi_{xy}(y)$ in $U_x$ is $\lambda_y(H_y)$. We can define:

$$
\hat N=\{ (x,(H)), H\subset H_x, (H)=(H^y_x)\}
$$

 where $x\in N$, $(U_x,H_x,\phi_x)$ is a local chart at $x$.
  We denote by $H"_x$ the set of subgroups of $H_x$ which are type
of an orbit, and by $H'_x$ the set of conjugacy classes of these
subgroups. Remark that the argument in Audin [4] p. 17 proposition
2.2.3 implies that we can assume that the number of types of
orbits contained in every chart is finite. The reunion $D_H$ of
the orbits whose type is $(H)$ is a submanifold. The following
proposition is shown for orbifolds by Chen and Ruan [11] p.7.

\medskip

{\bf Proposition 5.2.1.}

{\it There exists a generalized orbifold structure on $\hat N$.
Let $(U_x,H_x,p_{H_x})$ be a chart of $N$, and $(H)\in H'_x$. We
denote by $U^H$ the fixed point subset of $U$ by the action of
$H$, and by $C(H)$ the normalizer of $H$ in $H_x$, then
$((U^H,C(H)),C(H),\phi_H)$ is a chart of the generalized orbifold
$\hat N$, where $\phi_H:U^H\rightarrow U^H/C(H)$ is the natural
projection.}

\medskip

{\bf Proof.}
 Consider $(U_x,H_x,\phi_x)$ a chart at $x$. Let $y$
be an element of $\phi_x(U_x)$. Consider a chart
$(U_y,H_y,\phi_y)$, such that $U_y$ contains  an element $y'$ such
that $\phi_y(y')=y$ and $H_y(y')=y'$. The equivariant transition
function $\phi_{xy}:p_y(U_x\times_NU_y)\rightarrow
p_x(U_x\times_NU_y)$ where $p_x:U_x\times_NU_y\rightarrow U_x$ is
the canonical projection induces a morphism
$\lambda_{y'}:H_y\rightarrow H_x$. Let $H=H_{y'}^z$ and $h\in H$,
the element $\lambda_{y'}(h)$ fixes $\phi_{xy}(z)$.  We deduce a
map $\Phi$ which associates to $(y,(H))$ the projection of
$\phi_{xy}(y')$ in $\bigcup_{H=H_x^y\in H"_x}U^H_x/H_x$, where an
element $h$ of $H_x$ acts on $\bigcup_{H=H_x^y\in H"_x}U^H_x$ by
sending the element $c\in U_x^H$ to $h(c)\in U_x^{hHh^{-1}}$.

If instead of taking $H$, we take the element $H'$  in $(H)$,
$H'=aHa^{-1}$, $\phi_{xy}(az)\in U^{\lambda_{y'}(H')}$, and
$\Phi(y,(aHa^{-1}))$ is the projection to $\bigcup_{H=H_x^y\in
H"_x}U^H_x/H_x$ of $\phi_{xy}(y')$ in
$U_x^{\lambda_{y'}(hHh^{-1})}$.

If we take $y"$ such that $\phi_x(y')=\phi_x(y")$, $y"=by', b\in
H_x$, and $\phi_{xy}(bz)\in U^{b\lambda_{y'}(H)b^{-1}}$, and
$\Phi(y,(H))$ is the projection of $y"\in
U^{b\lambda_{y'}(H)b^{-1}}$ to $\bigcup_{H=H_x^y\in
H"_x}U^H_x/H_x$. Thus the map $\Phi$ is well defined. This map is
surjective; this can be shown by the fact that we can linearize
the action of compact Lie group. It is injective: If
$\phi(y,(H))=\phi(y_1,(H_1))$, and $\Phi(y,(H))$ and
$\Phi(y_1,(H_1))$ are the projections of $y'$ and $y_1'$ in
$\bigcup_{H=H_x^y\in H"_x}U^H_x/H_x$, there exists $a\in H_x$ such
that $y_1'=ay'$. This implies that $y=y_1$. The definition of
$\Phi$ implies then that $(H)=(H')$. Remark that the image of the
previous map is in bijection with $\bigcup_{(H)\in H'_x}
U^H/(C(H))$.

We endow $\hat N$ with the topology the topology generated by the
image of the maps $U^H\rightarrow \hat N$. The triples
$(U^H,C(H),\phi_H)$ defines a covering atlas of the generalized
orbifold where $\phi_H:U^H\rightarrow U^H/C(H)$ is the projection
map $\bullet$

\medskip

Let $H=H_x^y$, $U^H/H$ is an open subset of a suborbifold of $\hat
N$ completely determined by $(H)$ if $N$ is connected that we
denote $N_H$.

\medskip

Consider a pseudo-complex structure defined on $C_N$, this is
equivalent to suppose that each chart is endowed with a
pseudo-complex structure, and morphisms in $C_N$ preserve
pseudo-complex structures. Consider a chart $(U_x,H_x,\phi_x)$.
For every and $(H)$ in $H'_x$, we define
$2i_{(H)}=dim_{{\C}}(U_x)-dim_{{\C}}(N_H)$.

We can define:

$$
H^d(N)=\oplus H^{d-2i_{(H)}}(N_{H}).
$$

\section{ Sheaf of categories and gerbes in differentiable
categories.}

Recall that if  $C$ is a differentiable category endowed with a
topology, $U$  an object of $C$ and $R$ a sieve in $J(U)$. The
forgetful functor from $R$ to $C$ which sends a map $V\rightarrow
U$ to $V$ is Cartesian.

\medskip

{\bf Definition 6.1.}

Let $F:P\rightarrow C$ be a differentiable fibered functor, where
the category $C$ is equipped with a Grothendieck topology, we say
that $F$ is a sheaf of categories, if for every object $U$ of $C$,
and for every sieve $R\in J(U)$, the natural restriction map:

$$
Cart_C(C\uparrow U,F)\rightarrow Cart_C(R,F)
$$
is a $2$-descent map, otherwise said, an equivalence of
categories. (See Giraud [14])

\medskip

The sheaf of categories is called a gerbe bounded by the sheaf $H$
if the following conditions are satisfied:

$F$ is locally connected: this is equivalent to saying that for
every object $U$ of $C$, there exists a sieve $R\in J(U)$ such
that for every map $V\rightarrow U\in R$, the objects of the fiber
$P_V$ of $V$ are isomorphic each other.

There exists a sheaf in groups $H$ defined on $(C,J)$ such that
for every object $U\in C$, and $e_U\in P_U$ the group $Aut_U(e_U)$
of automorphisms of $e_U$ over the identity of $U$ is isomorphic
to $H(U)$, and these family of isomorphisms commute with morphisms
between objects and restrictions. The sheaf $H$ is called the band
of the gerbe.

\medskip

Let $(C,J)$ be a site, two fibered categories
$F_i,i=1,2:C_i\rightarrow C$ are equivalent, if there exists a
Cartesian isomorphism between $C_1$ and $C_2$.

An equivalence between the gerbes $F_i,i=1,2:P_i\rightarrow C$ is
a Cartesian isomorphism which commutes with their bands. Let $H$
be a sheaf defined on the differentiable site $(C,J)$, we denote
by $H^2(C,H)$ the set of equivalences classes of $H$-gerbes. This
set is often called the non-abelian 2-cohomology group of the
sheaf $H$.

\subsection{ The classifying cocycle of a gerbe.}

Suppose that the differentiable category $C$  has inductive
limits, finite projective limits, a final and initial object.

Let $R$ be a  covering of the final object $e$. We suppose that
$R$ is a good covering, that is every gerbe defined on an object
$X_i$  of $C$ such that there exists a map $X_i\rightarrow e$ in
$R$ is trivial and connected.

Let $F:P\rightarrow C$ be a gerbe, and $e_i$ an object of the
fiber $P_{X_i}$. There exists an isomorphism:

$$
u_{ij}:e^i_j\rightarrow e^j_i
$$

We denote by $c_{ijl}$ the isomorphism $u^j_{li}\circ
u^l_{ij}\circ u^i_{jl}$.

We have the relation:

$$
c^{i_2}_{i_1i_3i_4}u^{i_1i_2}_{i_4i_3}c^{i_4}_{i_1i_2i_3}u^{i_1i_2}_{i_3i_4}=c^{i_3}_{i_1i_2i_4}c^{i_1}_{i_2i_3i_4}
$$

The family of $2$-chains $c_{i_1i_2i_3}$ which satisfies the
relation above is called a non-abelian $2$-cocycle. Giraud [14]
has shown that there exists a 1 to 1 correspondence between the
set of gerbes bounded by $H$ and non abelian $H$ $2$-cocycles (see
also the proof in Brylinski [8]  p. 200-203 for commutative
gerbes).

\section{\bf Examples of sheaf of categories and gerbes.}

The differentiable category $C_H$ which represents the geometric
obstruction to extend the structural group of a principal bundle
is a gerbe (see page 17).

\medskip

Recently, Lupercio and Uribe [21] have introduced Abelian gerbes
on orbifolds. For an orbifold $N$, we can define a gerbe on $N$ to
be a gerbe defined on the Grothendieck site $J_N$ ( see p. 20).

\subsection{ Gerbes and $G$-structures.}

We are going to define a fundamental example of a gerbe, that we
are going to apply to the study of $G$-structures. Let $G$ be a
Lie group, and $H$ a closed subgroup of $G$. Consider a principal
$G$-bundle $p:P\rightarrow N$ over the manifold $N$. A natural
question is to ask wether the bundle has an $H$-reduction, that is
wether there exists coordinates change which take their values in
$H$. This problem is equivalent to the following question:
Consider the bundle $p':P'\rightarrow N$ whose typical fiber is
the homogeneous space $G/H$ obtained by making the quotient of
each fiber of $p$ by $H$. Is there exists a global section of
$p'$?(see Albert and Molino [2] p. 64). We have the following:

\medskip

{\bf Proposition 7.1.1.}

{\it The correspondence defined on the category of open subsets of
$N$, which assigns to every open subset $U$ the category $C_H(U)$,
whose objects are $H$-reductions of the restriction of $p$ to $U$,
and whose morphisms, are morphisms of $H$-bundles is a sheaf of
categories.}

\medskip

{\bf Proof.}
 Gluing condition for objects.

Let $(U_i)_{i\in I}$ be an open covering of $U$, $e_i$, an object
of $C_H(U_i)$ such that there exists a morphism
$u_{ij}:e^i_j\rightarrow e^j_i$ such that
$u_{ij}^lu_{jl}^i=u_{il}^j$. Then there exists an object $e$ in
$C_H(U)$ whose restriction to $U_i$ is $e_i$, since we can glue
$H$-bundles.

Gluing conditions for arrows:

Let $e$ and $e'$ be two objects of $C_H(U)$, the correspondence
which assigns to every open subset $V$ of $U$,
$Hom_{C_H(U)}(e_{\mid V},e'_{\mid V})$ is a sheaf, since it is the
sheaf of morphisms between two bundles $\bullet$

\medskip

This sheaf of categories can be applied to the following
situation: suppose that $N$ is a $n$-dimensional manifold. Let
$R_p(N)$ be the bundle of $p$-linear frames of $N$, and $G$ a
subgroup of $Gl_p(n,{\R})$ the group of invertible $p$-jets of
${\R}^n$. The geometric obstruction of the existence of a
$G$-structure on $N$ is defined by the sheaf of categories that we
have just defined.

Let $U$ be a contractible open subset of $N$, $C_G(U)$ is not
empty, since the restriction of $P$ to $U$ is a trivial bundle.
But the objects of $C_G(U)$ are not always isomorphic: suppose
that $N$ is a $n$-dimensional manifold, and take $G=O(n,{\R})$;
the $G$-reductions of the bundle of linear frames $R(N)$ of $N$
define the differentiable metrics. It is  well-known  that two
differentiable metrics are not locally isomorphic if their
curvatures are distinct.

A particular situation is the example of flat $G$-structures like
symplectic structures (see Albert and Molino [2] p. 177). For
every elements $x$ and $y$ in $N$, there exists neighborhoods
$U_x$ and $U_y$ of $x$ and $y$ in $N$, and a diffeomorphism
$h:U_x\rightarrow U_y$ which preserves the $G$-structures induced
by $N$ on $U_x$ and $U_y$.  If $U$ is contractible open subset of
$N$, two elements of $C_G(U)$ are connected. The theory of gerbes
and $G$-structures, will be intensively studied in [42].

\subsection{ Gerbes and invariant scalar product on Lie groups.}

Another example of gerbes can be described as follows: consider a
Lie group $H$ which is not commutative, and $L$ a lattice in $H$.
Consider  the manifold $H/L$, and let suppose that $H$ is endowed
with an orthogonal bi-invariant metric: this is equivalent to the
existence of a scalar product $<,>$ (i.e a non-degenerated real
valued bilinear form not necessarily positive definite) on the Lie
algebra ${\cal H}$ of $H$ such that for every elements $x,y,z\in
{\cal H}$:

$$
<[x,y],z>+<y,[x,z]>=0.
$$

 The $3$-invariant form $\nu$ defined on the Lie algebra ${\cal
H}$ of $H$ by:

$$
\nu(x,y,z)=<[x,y],z>
$$
 defines on $H/L$ a closed  $3$-form $\nu_L$.

 \medskip

 The space of bilinear symmetric forms on ${\cal H}$ corresponds to real
 $3$-cocycles as shows Koszul [17] p. 95. Medina [26] has shown that
 the dimension of this space is either $1$ or $2$.

\medskip

Let $H$ be a $n$-dimensional nilpotent Lie group, and $L$ a
lattice of $H$. Recall that there exists a basis $e_1,...,e_n$ of
the Lie algebra ${\cal H}$ of $H$, such
$[e_i,e_j]=\sum_{ijl}c_{ijl}e_l, c_{ijl}\in {\Q}$ (see Raghunathan
[35] p. 34). We say in this situation that the constants of
structure are rational. A lattice $L$ is the image of
${\Z}e_1\oplus...\oplus{\Z}e_n$ by the exponential map.

\medskip

{\bf Proposition 7.2.1.}

{\it Under the notations above, if there exists an invariant
scalar product $<,>$, such that $<e_i,e_j>\in {\Q}$, then  the
$3$-form  $\nu$ on $H/L$ induces canonically a rational $3$-form
$\nu_L$ on $H/L$.}

\medskip

{\bf Proof.}
 The proof uses a theorem of Nomizu quoted in
 [35] Raghunathan p. 123 in the real case. Let $H_0$ be the center of
 $H$, the intersection $L_0=L\cap H_0$ is a lattice in $H_0$. (If $H$ is commutative,
 we take $H_0$ to be a non trivial subgroup different of $H$.
 See  Raghunathan [35] p. 40). Thus the foliation of $H/L$ whose leaves are
 orbits of $H_0$, has compact leaves. The space of leaves of this foliation is
 the quotient $M$ of $H/H_0$ by $L/L_0$. The natural projection $H/L\rightarrow M$
  is a fibration whose fibers are $n$-dimensional torus $T^n$, where
  $n$ is the dimension of $H_0$.

  We can apply the Leray-Serre spectral sequence to this fibration
  for the rational cohomology we obtain:

  $$
  E_2^{p,q}=H^p(M,H^q(T^n,{\Q}))\simeq H^p(M,\Lambda{\Q}^q),
  $$

  $$
  E^{p,q}_{\infty}\Longrightarrow H^{p+q}(N,{\Q})
  $$

  Consider ${\hat E}^{p,q}_*$ the Leray-Serre spectral sequence associated
  to the space of $H$-invariant forms on $N$ and $M$, we have:

  $$
  {\hat E}_2^{p,q}=H^p({\cal H}/{\cal H}_0, \Lambda{\Q}^q),
  $$

  $$
  \hat{E}^{p,q}_{\infty}\Longrightarrow H^{p+q}({\cal H},{\Q}).
  $$

  The recursive hypothesis implies that $H^*(M,{\Q})=H^*({\cal H}/{\cal H}_0,{\Q})$.
This implies that $H^*(N,{\Q})=H^*({\cal H},{\Q})$. The image of
$\nu$ by the isomorphism $H^3({\cal H},{\Q})\rightarrow
H^3(N,{\Q})$ is the form $\nu_L$. The result of Koszul [17] p. 95
shows that we can realize this form by using an invariant bilinear
form $\bullet$

\medskip

 The classification theorem of Giraud [14]  implies the existence of a gerbe over
 $H/L$ whose classifying class is the cohomology class of $p\nu_L$, where $p$ is an integer.
 We call such a gerbe, a Medina-Revoy gerbe.

 \bigskip

 {\bf Examples of Medina-Revoy gerbes.}

\medskip

 Lie groups endowed with bi-invariant scalar product have been intensively
 studied by Aubert, Dardie, Diatta, Medina and Revoy. Medina and
 Revoy [27] have shown that they can be constructed from simple Lie groups
 and the $1$-dimensional Lie group by the processus of double
 extension. Here is an example of a Medina Revoy gerbe constructed
 from the double extension of the two dimensional Euclidean space,
 endowed with its commutative structure of a Lie algebra.

 Consider the nilpotent Lie algebra constructed as follows: Let
$(e_1,e_2)$ be an orthogonal basis of the $2$-dimensional
Euclidean space $(U,<,>)$, and $h:U\rightarrow U$ the linear
endomorphism such that $h(e_1)=e_2, h(e_2)=0$ considered also as a
derivation of the trivial underlying Lie algebra of $U$. Let $V$
be the $1$-dimensional commutative Lie algebra, and $V^*$ its
dual. For every elements $u_1,u_2\in U$, we denote by
$w(u_1,u_2):V\rightarrow V^*$ the map which assigns to $v\in
V={\R}$ the scalar $<vh(u_1),u_2>$. The double extension of
$(U,<,>)$ by $V$ and $h$ is the nilpotent Lie algebra ${\cal
L}=V^*\oplus U\oplus V$ whose bracket is defined by the formula:

$$
[(v'_1,u_1,v_1);(v'_2,u_2,v_2)]=(w(u_1,u_2),v_1h(u_2)-v_2h(u_1),0)
$$

The Lie algebra $V^*\oplus U\oplus V$ is endowed with the scalar
product:

$$
<(v'_1,u_1,v_1);(v'_2,u_2,v_2)>'=<u_1,u_2>+v_1v_2+v'_1(v_2)+v'_2(v_1)
$$

The constant of structures of ${\cal L}$ are integral in its
canonical basis. The $3$-form $\nu_L$ defined on $L$ by
$(u,v,w)\rightarrow <[u,v],w>'$ is rational.

 Let $\Gamma$ be the lattice of the $1$-connected Lie group $L$
associated to ${\cal L}$ which is generated by the image of an
integral basis of ${\cal L}$. The classification theorem of Giraud
implies the existence of a gerbe on $L/\Gamma$ whose classifying
cocycle is $p\nu_L$, where $p\in {\N}$ is such that $p\nu_L$ is
integral.

\subsection{\bf A sheaf of categories on an orbifold with
singularities.}

Let $H$ be a compact Lie group which acts on a manifold,  the
quotient space $N/H$  is an example of a generalized orbifold
$C_N$ (see definition 3.2.1).

\medskip

{\bf Proposition 7.3.1.}

{\it Let $N$ be a generalized compact orbifold. The correspondence
defined on the category of open subsets of $N$  which assigns to
$U$ the category $C_N(U)$, whose objects are elements
$(P,H_P,\phi_P)$ of $C_N$, such that the image of $\phi_P$ is $U$
is a sheaf of categories.}

\medskip

{\bf Proof.}
 Gluing conditions of objects.

Let $U$, be an open subset of $N$, and $(U_i)_{i\in I}$ an open
covering of $U$. Consider for each $i\in I$, an object
$e_i=(P_i,H_i,\phi_i)$ in $C_N(U_i)$, and a morpism
$u_{ij}:e^i_j\rightarrow e^j_i$ such that
$u_{ij}^lu_{jl}^i=u_{il}^j$. Since the morphisms $u_{ij}$ are
local diffeomorphisms, there exists a manifold $P$ obtained by
gluing the family of manifolds $P_i$ with $u_{ij}$. We can glue
the Lie groups $H_i$ and their actions to define a Lie group $H$
which acts on $P$, and such that the map $P/H\rightarrow N$ is a
local homeomorphism:

Let $l_i$ be the restriction  of the action of $H_i$ to
$p_i(P_i\times_NP_j)$. We can identify $p_i(P_i\times_NP_j)$ with
$p_j(P_i\times_NP_j)$ with $u_{ij}$. We denote $H_{ij}$ the limit
of the maps $l_i$ and $l_j$. The Lie group $H_{ij}$ acts on the
gluing of $P_i$ and $P_j$ by $u_{ij}$. Without restricting the
generality, we can suppose that $I$ is a numerable set, construct
$H_{01}$, $H_{01..n}$ obtained by gluing recursively the action of
$H_0,...,H_n$. The Lie group $H$ is the limit of the groups
$H_{01..n}$.

Gluing condition of arrows.

Let $U$ be an open subset of $N$, and $P$, and $P'$ two objects of
$C_N(U)$. The correspondence defined on the category of open
subsets of $U$, which assigns to $V$ the set $Hom_{C_N(V)}(P_{\mid
V},P'_{\mid V})$, where $P_{\mid V}$ is $\phi_P^{-1}(V)$ is a
sheaf since we can glue differentiable maps $\bullet$

\medskip

Another example of sheaf of categories is defined by the theory of
foliages (see Molino [31], or definition 3.3.1). Let $N$ be a
topological manifold  endowed with a structure of a foliage. A
natural problem is to determine wether this foliage is induced by
a foliation on a manifold. The following proposition provides the
obstruction which solves this problem.

\medskip

{\bf Proposition 7.3.2.}

{\it Let $N$, be a topological space endowed with the structure of
a foliage, for every open subset $U$ of $N$, denote by $C_N(U)$
the class of objects $(U, V, {\cal F},\pi)$ of $C_N$. The
correspondence which assigns $C_N(U)$ to $U$ is a sheaf of
categories which is the geometric obstruction of the existence of
a manifold $\hat N$, endowed with a foliation ${\cal F}_N$, such
that $N$ is the space of leaves of ${\cal F}_N$.}

\medskip

{\bf Proof.}
 Gluing conditions of objects.

Let $U$ be an open subset of $N$, $(U_i)_{i\in I}$ an open
covering of $U$ such that for each element $i$ of $I$, there
exists an object $e_i=(U_i,V_i,{\cal F}_i,\pi_i)$ in $C_N(U_i)$,
morphisms $u_{ij}:e^i_j\rightarrow e^j_i$ such that
$u_{ij}^lu_{jl}^i=u_{il}^j$. The morphisms $u_{ij}$ allow to glue
the family of manifolds $V_i$ to obtain a manifold $V$, on which
is defined a foliation ${\cal F}$ whose restriction to $U_i$ is
${\cal F}_i$.

Gluing condition for arrows.

Let $e=(U,V,{\cal F},\pi)$ and $e'=(U,V',{\cal F}',\pi')$ two
objects of $C_N(U)$. The correspondence defined on the category of
open subsets of $U$, which assigns to $U'$ the set
$Hom_{C_N}(e_{\mid U},e'_{\mid U})$, is a sheaf, since we can glue
differentiable foliated maps $\bullet$

\medskip

There exists a sheaf $L$ on $N$ which assigns to every open subset
$U$ of $N$, the set of isomorphisms of an object $e_U$ of $C_N(U)$
(the foliated isomorphisms), and the sheaf of categories that we
have just defined is a gerbe bounded by $L$.

\section{\bf Differential geometry of sheaves of categories.}

In this part, we are going to analyze the tools defined in the
general context of differentiable categories to study their
geometry by using the underlying topology.

\medskip

Let $C$ be a differentiable category endowed with a topology. We
suppose that  $C$ has a final object and a good global covering
$(U_i)_{i\in I}$ (see p. 24). Let $P\rightarrow C$ be a gerbe, we
suppose that there exists a Lie group $H$,  a principal $H$-torsor
$A:Aut(P)\rightarrow C$, such that  every object $e_U\in P_U, U\in
C$ is a bundle $p_{e_U}:e_U\rightarrow U$ endowed with a free
right action of $H$. The set of morphisms between two objects of
$P_U$ are morphisms between bundles which project to the identity
of $U$, and the set of automorphisms of $e_U$ can be identified
with gauge transformations of $Aut(P)(U)$ by a map which commutes
with morphisms between objects and with restrictions. We denote by
$aut(P)$ the vector ${\cal H}$-bundle on $C$ associated to
$Aut(P)$: If the coordinate changes of $Aut(P)$ are defined by the
maps $(u_{ij})_{i,j\in I}$, the coordinate changes of $aut(P)$ are
defined by the map $(Ad(u_{ij}))_{i,j\in I}$. Such a gerbe is
called a $H$-gerbe.

\subsection{ Induced gerbes.}

Let $p:P\rightarrow C$ be an $H$-principal gerbe, that is: for
every object $U$ of $C$, the map $p_{e_U}:e_U\rightarrow U$,
endows $e_U$ with the structure of a $H$-principal bundle.
Consider a morphism of Lie groups $h:H\rightarrow H'$, we can
construct a principal $H'$-gerbe $p':P'\rightarrow C$ as follows:

Let $a:A\rightarrow U, U\in C$ be an object of $P_U$, it is a
$H$-principal torsor defined by the trivialization $(U_i,u_{ij}\in
H)_{i,j\in I}$. We can define the image of $a$ by $h$. It is the
torsor whose coordinates change are defined by:
$(U_i,h(u_{ij}))_{i,j\in I}$. The family of these images is the
induced gerbe.

\medskip

An example is the situation when $H$ is $SU(n)$ or $O(n)$, and $h$
is the determinant morphism.

\medskip

{\bf Proposition 8.1.}

{\it There exists a connective structure on each  $H$-gerbe
$p:P\rightarrow C_N$, where $N$ is a manifold, and $C_N$ the
differentiable category associated to $N$ defined at p.5.}

\medskip

{\bf Proof.}
 Let $(U_i)_{i\in I}$ be an open cover of $N$, and
$e_i$ an object of $P_{U_i}$, we denote by $e'_i$ the quotient of
$e_i$ by the action of $H$.

Consider an isomorphism $u_{ij}: e^i_j\rightarrow e^j_i$. We can
construct from $u_{ij}$  a morphism $u'_{ij}:{e'}^i_j\rightarrow
{e'}^j_i$ which is its quotient. With these morphisms, we can glue
the family  of quotients $(e'_i)_{i\in I}$ to define a fiber
bundle $p':P'\rightarrow N$.

Consider a connection on $p'$: this is a distribution ${\cal D}'$
of $P'$ whose rank is the dimension of $N$, and which is
transverse to the fiber of $p'$ (compare with McDuff and Salamon
[23] p. 210). The distribution ${\cal D}'$ can be also defined by
a $1$-form $\theta$ on $P'$ which takes its values in the
canonical bundle over $P'$, such that for each $x\in P'$, the
fiber of this bundle is the tangent space of $P'_x$, the fiber at
$p'(x)$. We suppose that if $v\in T_xP'$, $x\in U$, $\theta(v)=v$.
Such a distribution can be constructed by using a differentiable
metric on $P'$, and by taking the orthogonal of the fiber.

Let $U$ be an open subset of $N$, and $e_U$ an object of $P_U$. We
denote by $Co(e_U)$ the set of  $1$-forms defined on the  bundle
$e_U\rightarrow U$ which take their values in the canonical
vector bundle over $U$ whose fiber at $x$ is the tangent space of
the fiber of $e_U\rightarrow U$ at $x$. We suppose  that for every
$\alpha\in Co(e_U)$, and every element $A\in {\cal H}$ which
generates the fundamental vector field $A^*$, we have:

$$
\alpha(A^*)=A.
$$

We suppose also that the each element of $Co(e_U)$  descends to
the restriction of $\theta$ to $U$.

 Such a connection can be
constructed as follows: Let $(U_i)_{i\in I}$ be a trivialization
of $e_U\rightarrow e_U/H$. We can define on $U_i\times H$ a
distribution invariant by $H$ whose projection to $U_i$ is the
restriction of ${\cal D}'$. Such a distribution is defined by a
$1$-form $\alpha_i$ whose projects to the restriction of $\theta$
to $U_i$. By using a partition of unity, we deduce that $Co(e_U)$
is not empty $\bullet$

\subsection{ Reduction to a situation similar to the motivating example.}

Now we reduce the study of the differential geometry of a
$H$-gerbe to a situation similar to the motivating example (see p.
17) by using the following construction. Consider a $H$-gerbe
$p:P\rightarrow N$. We have seen that there exists a fiber bundle
$p':P'\rightarrow N$ such that for a good covering $(U_i)_{i\in
I}$ of $N$, the restriction of $P'$ to $U_i$ is the quotient of
every object $e_i$ of $P_{U_i}$ by $H$. Let $F$ be the fiber of
this bundle. We can define the bundle $L(p'):L(F)\rightarrow N$
such that for every element $u\in N$ the fiber of $L(F)$ at  $u$
is the set of linear frames of its fiber, $F_u$. Let $U$ be an
open subset of $U$, and $e_U$ an object of $P_U$. We can define
the pullback of the maps $L(F)_{\mid U}\rightarrow e_U/H$ and
$e_U\rightarrow e_U/H$ which is an $H$-principal bundle $e'_U$
over the restriction $L(F)_{\mid U}$. The class $P'_U$ whose
elements are the $e'_U$ just constructed defines a gerbe
$L(P)\rightarrow N$. This allows to deal only with to a situation
similar to the motivating example in supposing that the bundle
$P'\rightarrow N$ is principal. In fact in the sequel, we deal
only with the motivating example. Moreover we suppose that
connective structures are constructed by using connections of $P'$
like at the proposition 8.2.1.

Another variant of the previous construction is the following: let
$p:P\rightarrow N$ be a $H$-principal bundle defined on the
manifold $N$, to each cocycle $c\in H^2_{Cech}(N,H)$ we can
associate a gerbe $C$ whose classifying cocycle is $c$ (see Giraud
[14]). The previous construction allows to construct a gerbe
$R(C)$ above the frames bundle of $N$, such that an object of
$R(C)(U)$ is the pull back of an object of $C(U)$ by the
projection map $R(U)\rightarrow U$, where $R(U)$ is the bundle of
linear frames of $U$. Thus we can define connective structures
above connections of $R(N)$.

\medskip

{\bf Definition 8.2.1.}

A curving (see also Brylinski [8] p. 211) defined on the
connective structure $Co$ of the $H$-gerbe $p:P\rightarrow C$, is
a map which assigns to every object $e_U$ of $P_U$, and every
element $\theta\in Co(e_U)$, a $2$-form  $L(e_U,\theta)$ which
takes its values in ${\cal H}$ such that the following properties
are satisfied:

Let $h:e_U\rightarrow e_{U'}$ be a morphism, for every $\theta\in
Co(e_{U'})$, we have:

$$
L((e_{U}),h^*(\theta))=h^*(L(e_{U'},\theta)).
$$

If $h$ is an automorphism of $e_U$, and $\theta$ an element of
$Co(e_U)$, we have:

$$
L(e_U,h^*(\theta))=Ad(h^{-1})(L(e_U,\theta))
$$

 Let $\alpha$ be an element of $\Omega^1(U,aut(U))$, we have:

$$
L(e_U,\theta+\alpha)=L(e_U,\theta)+d(\alpha)+{1\over
2}([\alpha,\alpha]+[\theta,\alpha]+[\alpha,\theta]).
$$

The correspondence $(e_U,\theta)\rightarrow L(e_U,\theta)$ is
natural in respect to restrictions and morphisms between objects.

\medskip

Remark that since the gerbe considered here is associated to the
motivating example, the $2$-form
$L(e_U,\theta+\alpha)-L(e_U,\theta)$ is ${\cal H}$-valued.

\bigskip

{\bf Proposition 8.2.1.}

{\it Let $Co$ be a connective structure defined on the
$H$-principal gerbe $p:P\rightarrow C$, there exists a curving.}

\medskip

{\bf Proof.}
 Compare the following proof with Brylyinski [8] p.
212). We are going to assume that there exists a $L'$-bundle
$p':P'\rightarrow N$, an exact sequence of Lie groups
$1\rightarrow H\rightarrow L\rightarrow L'\rightarrow 1$ such that
the gerbe is the geometric obstruction to lift the structural
group $L'$ of $p'$ to $L$. We suppose also that the connective
structure defined on the principal gerbe is constructed as in the
proposition 8.1. Thus there exists a connection $\alpha$ defined
on the principal bundle $p':P'=P/H\rightarrow N$ such that for
every open subset $U$  of $N$, $e_U$ an object of $P_U$, the
elements of $Co(e_U)$ are connections which projects to $\alpha$.

Let ${\cal H}$, ${\cal L}$, and ${\cal L}'$ be the respective Lie
algebras of $H$, $L$ and $L'$. Let $u$ be a linear section of the
canonical map ${\cal L}\rightarrow {\cal L}'$. We can define the
form $\alpha_0=u\circ\alpha$ on $P'$. Let $\theta$ be an element
of $Co(e_U)$, consider the form $\alpha_{e_U}$, the pull-back of
the restriction of $\alpha_0$ by the canonical projection
$e_U\rightarrow e_U/H=P'_{\mid U}$. We set:

$$
L(e_U,\theta)= d\theta +{1\over
2}[\theta,\theta]-(d\alpha_{e_U}+{1\over
2}[\alpha_{e_U},\alpha_{e_U}]).
$$
$\bullet$

\bigskip

{\bf Definition 8.2.2.}

Let $L$ be a curving of the connective structure $Co$ defined on
the $H$-gerbe $P\rightarrow C_N$, where $N$ is a manifold, and
$C_N$ the canonical differentiable category associated to $N$ (see
page. 5). Let $(U_i)_{i\in I}$ be a good covering of $N$, and
$e_i$ and object of $P_{U_i}$, and $\alpha_i$ an element of
$Co(e_i)$ and $u_{ij}:e^i_j\rightarrow e^j_i$ a morphism. Let
$\Omega_i$ be the curvature of $\alpha_i$. On $U_i\cap U_j$, the
difference $\Omega^i_j-u_{ij}^*\Omega^j_i$ defines a $1$-Cech
$aut(P)$-cocycle. The  cohomology class does not depend of the
elements $e_i$ in $P_{U_i}$ and of the elements $\alpha_i$ used to
define it, since $Co(e_i)$ is a torsor whose vector space is a
space of ${\cal H}$-valued $1$-forms.

 The DeRham-Cech isomorphism allows to identify this cocycle
to a $3-aut(P)$ form $D$ called a curvature of the connective
structure.

\bigskip

{\bf Definition 8.2.3: Holonomy.}

\medskip

We are going to reduce this definition to the commutative case,
all the groups are compact and complex, as well as the vector
bundles. Let $p:P\rightarrow C$ be a differentiable $H$-gerbe,
endowed with a connective structure and a curving $L$. Suppose
that $P$ is the geometric obstruction to lift a $G$-bundle $P'$
over $N$ to a $G'$-bundle given the exact sequence of Lie groups
$1\rightarrow H\rightarrow G'\rightarrow G\rightarrow 1$. We
assume also that in fact, the $G$-bundle $P'$ is the frames bundle
of a vector bundle $V$ over $N$, and the objects of $P$ are also
associated to vector bundles.

Let $N$ be a surface, Consider a morphism $h:N\rightarrow C$. We
can pull-back $p$ by $h$ and obtain a  gerbe $p_N:P_N\rightarrow
N$, endowed with a connective structure which has a curving. The
quotient of $p_N$ by $H$ is the pull-back $p'_G:P^G_N\rightarrow
N$ of $P'$ by $h$, which is a reduction of the frames bundle of
the pull-back $V_N$ of $V$ by $h$. A well-known result implies
that $V_N$ is isomorphic to the summand of complex line bundles
(see McDuff and Salamon [23] p. 80). Thus we can assume that the
structural group of $P_N$ is a commutative subgroup $G_N$ of $G$.
The gerbe $p_N$ is also the geometric obstruction to extend the
structural group $G_N$ to $G'_N$ given an exact sequence of Lie
groups $1\rightarrow H_N\rightarrow G'_N\rightarrow G_N\rightarrow
1$. This implies that  we can assume also that $H_N$ and $G'_N$
are commutative.

 The pull-back of the connective structure and the curving of $p$,
 induces a connective structure $Co_N$ of $p_N$ and a curving
 $L_N$, moreover the connection $\alpha$ used to construct the
 connective structure of $p$ (see proposition 8.2.1) is supposed to be Hermitian, as well as the elements of
 $Co(e_U)$, where $e_U$ an object of the gerbe. Thus the connection
   $\alpha_N$ on $P^{G_N}_N$, which induces the connective
 structure $Co_N$ preserves every Hermitian reduction. We just have to recall the definition of the
 holonomy in the commutative case.

Let $(U_i)_{i\in I}$ be a good cover of $N$. Let $e_i$ be an
object of ${P_N}_{U_i}$, and let $\theta_i$ be an element of
$Co_N(e_i)$, we can suppose that this connection takes its values
in the Lie algebra of $G'_N$. We denote by $L(e_i,\theta_i)$ the
curving associated to the element $\theta_i\in Co(e_i)$. Denote by
$\theta_{ij}$ the form $\theta^i_j-u_{ij}^*(\theta^j_i)$. Since
$N$ is $2$-dimensional, there exists a $1$-form $h_i$ such that
$dh_i=L(e_i,\theta_i)$. We have:

$$
\theta_{ij}=h_j-h_i+da_{ij}
$$

We can set

$$d_{ijl}=c_{ijl}^{-1}a_{jl}^{-1}a_{il}a_{ij}^{-1}
$$

where $c_{ijl}$ is the classifying cocycle of $p_N$. The chain
$d_{ijl}$ is the holonomy cocycle of the gerbe $p_N$ endowed with
its connective structure.

The Cech-DeRham isomorphism allows to identifies this form with a
$2$-form $\Omega$ on $N$.
 (See also Mackaay and Picken [22] p. 27).

\subsection{Holonomy and functor on loops space.}

Consider the category $C^2$ whose objects are  maps:
$h:C^1+..+C^1\rightarrow N$, where $C^1+...+C^1$ is a finite
disjoint union of circles. A morphism between the objects $h$ and
$h'$, is a map from a surface $l:N\rightarrow C$ such that the
restriction of $l$ to the boundary of $N$ is the sum of the maps
$h$ and $h'$. The
 holonomy defines a functor $D$ on $C^2$ which associates to $h$
the complex line ${\C}$. Let $l:N\rightarrow C$ be a morphism
between $h$ and $h'$. The real holonomy around $N$ is the image
$D(l)$ of $l$ by $D$.

\bigskip

We suppose here that the structural group $H$ of the $H$-gerbe
$P\rightarrow N$ defined over the manifold $N$ and endowed with a
connective structure is contained in $Gl(n,{\C})$. We can relate
this gerbe which is associated determinant gerbe as follows:

\medskip

{\bf Proposition 8.3.1.}

{\it Suppose that $H$ is included in  $Gl(n,{\C})$, then  the
trace of the curvature of a principal $H$-gerbe $P\rightarrow N$
is the curvature of the associated determinant gerbe
$det(P)\rightarrow N$; that is the gerbe induced by the
determinant morphism $H\rightarrow {\C}$. (See 8.1).}

\medskip

{\bf Proof.}
 Let $p:P'\rightarrow N$ be the quotient of the gerbe
by $H$. For every open subset $U$ of $N$,  the objects of $P_U$
are $H$-principal bundles over $P'_{\mid U}$ the restriction of
$P'$ to $U$.

Let $det_{e_U}:e_U\rightarrow det(e_U)$ the determinant morphism
which associates to the object $e_U$ of $P_U$, the corresponding
object $det(e_U)$ in $det(P)$, and $\alpha\in Co(e_U)$ whose
kernel defines the distribution ${\cal C}_{\alpha}$ on $e_U$.
Since $det_{e_U}$ is an equivariant morphism between the
$H$-bundle $e_U\rightarrow P'_{\mid U}$ and $det(e_U)$, the image
of ${\cal C}_{\alpha}$ is a distribution $det({\cal C}_{\alpha})$
invariant by the action of $U(1)$ on $det(e_U)$ and transverse to
the fiber. Such a construction thus defines a connective structure
on the determinant gerbe.

 Suppose that the form $\alpha$ is ${\cal L}\oplus
gl(n,{\C})$-valued, where ${\cal L}$ is the Lie algebra of $L$ the
structural group of $P'$. Then the connection form which defines
the distribution $det({\cal C}_{\alpha})$ is the composition of
$\alpha$ and $(I_{{\cal L}},trace_{gl(n,{\C})})$. This implies
that the curving of $(det(e_U),{\cal C}_{e_U})$ is the trace of
the curving $L$. This implies the result $\bullet$

\subsection{Canonical relations associated to a connective
structure on a gerbe.}

\medskip

In this part we are going to determine canonical relations
associated to a connective structure. (compare with  Breen and
Messing [7] p. 58) The morphisms $u^*$ are pull-back, and the
morphisms the $u_*$ are inverse of pull-back. Let $e_i$ be an
object of $P_{U_i}$, and $\alpha_i\in Co(e_i)$. Consider the
restriction $e^i_j$ of $e_j$ to $U_{ij}$ and
$u_{ij}:e^i_j\rightarrow e^j_i$ an arrow. The $1$-form
${u_{ij}}_*(\alpha_j)$ is an element of $Co(e^j_i)$, since
$Co(e^j_i)$ is a torsor, there exists a $1$-form $\alpha_{ij}$
such that:

$$
\alpha^j_{i}={u_{ij}}_*(\alpha^i_j)+\alpha_{ij}
$$

We have seen that the family of forms $\alpha_{ij}$ verifies the
equations:

$$
{u_{ij}}_*\alpha_{jl}-\alpha_{il}+\alpha_{ij}=
{u_{il}}_*(\alpha^{ij}_l-Ad(c_{ijl}^{-1})(\alpha^{ij}_l))-c_{ijl}^{-1}dc_{ijl}
$$

 where $c_{ijl}$ is the map $u_{li}^ju_{ij}^lu_{jl}^i$.

Let $L$ be a curving of the connective structure. Denote by
$L_{ij}$ the $2$-form
$L_j(e^i_j,\alpha_j)-L_j(e^i_j,u_{ij}^*(\alpha_i))$. We have:

$$
L_{jl}-L_{il}+u^*_{jl}L_{jl}=
$$

$$
=L_l(e^{ij}_l,\alpha_l^{ij})-L_l(e^{ij}_l,{u_{jl}^i}^*(\alpha_j^{il}))
-(L_l(e^{ij}_l,\alpha_l^{ij})-L_l(e^{ij}_l,{u_{il}^j}^*(\alpha_i^{jl})))
+{u_{jl}^i}^*(L_j(e^{il}_j,\alpha_j^{il})-L_j(e^{il}_j,{u_{ij}^j}^*(\alpha_i^{jl})))
$$

$$
=L_l(e^{ij}_l,{u_{il}^j}^*(\alpha^{jl}_i))-L_l(e_l^{ij},{u_{jl}^j}^*{u_{ij}^l}^*(\alpha_i^{jl}))
$$

$$
={u^j_{il}}^*(L_i(e^{jl}_i,\alpha_i^{jl})-Ad({c'_{ijl}}^{-1})L(e^{jl}_i,\alpha_i^{jl}))
$$

where $c'_{ijl}=u_{ij}^lu_{jl}^iu_{li}^j$.

\subsection{\ Uniform distributions and gerbes.}

Another treatment of the differentiable structure on gerbes can be
done as follows: Let $p:P\rightarrow N$ be a $H$-gerbe defined on
a manifold. We reduce the study to the motivating example.  The
natural way to study the differential geometry of a principal
bundle is to use the theory of connections. Unfortunately,
connections defined on a principal bundle are not necessarily
invariant by the gauge group. This motivates the definition of a
torsor of connections, which is invariant by the automorphisms
group.

There exists another point of view used by Molino in his thesis
(see [29]). Molino has studied the notion of invariant
distributions on principal bundles. An invariant distribution on a
principal bundle  is a right invariant distribution. We do not
request here that the dimension of the distribution is the
dimension of the basis space of the bundle. The invariant
distribution is transitive, if its summand with the tangent space
of the fiber, generates the tangent space of the bundle, (see
Molino [29] p. 180), transitive distributions  are nothing but
equivalence classes of connections. We focus on the motivating
example of a $H$-gerbe $P\rightarrow N$; there exists an exact
sequence of Lie groups $1\rightarrow H\rightarrow L'\rightarrow
L\rightarrow 1$, a $L$-principal bundle $P'\rightarrow N$, such
that the gerbe $P\rightarrow N$ is the geometric obstruction to
lift the structural group of the previous principal bundle to
$L'$. The objects of $P_U$ are principal $H$-bundles over the
restriction of $P'$ to $U$. For a connection $\theta$ defined on
$P'$, we can define on each object $e_U$ of $P_U$, the transitive
distribution which is the kernel of the pull-back of $\theta_{\mid
P'_{\mid U}}$ to $e_U$. This transitive distribution is invariant
by the automorphisms of $e_U$.

\medskip

Suppose now that the extension which defines the lifting problem
is central. Let ${\cal H}$, ${\cal L}'$ and ${\cal L}$ be the
respective Lie algebras of $H$, $L'$ and $L$. The coordinate
changes $(u_{ij})_{i,j\in I}$ of the principal $L$-bundle
$P'\rightarrow N$, define a ${\cal L}$-bundle $P^{{\cal L}}$ over
$N$ whose coordinates changes are $(Ad(u_{ij}))_{i,j\in I}$. We
can lift $u_{ij}$ to an element $u'_{ij}$ of $L'$, since the
extension is central, we can define a ${\cal L}'$-bundle $P^{{\cal
L}'}$ over $N$ whose coordinates change are $(Ad(u'_{ij}))_{i,j\in
I}$. There exists a canonical projection $p_0:P^{{\cal
L}'}\rightarrow P^{{\cal L}}$.

 A connection structure defined on the gerbe $P$, is defined as follows:

 Let $U$ be an open subset of $N$, and $e_U$ an object of $P_U$, there exists
 a transitive distribution ${\cal D}_{e_U}$ of $e_U$ which is right invariant.

Let $h:e_U\rightarrow e_{U'}$ be a morphism in $P$, we assume that
the pull-back of ${\cal D}_{e_{U'}}$ by $h$ is ${\cal D}_{e_U}$.

  The distribution ${\cal D}_{e_U}$ is not assume to be uniform, see (Molino [29] p. 184)
  when it is uniform, the connection structure can be defined
 by a family of $1$-forms $\theta_{e_U}:e_U\rightarrow P^{{\cal L}'}$
 which verify the following conditions:

  if $x$ is an element of $e_U$, and $v$ an element of $T{e_U}_x$, the
tangent space of $e_U$ at $x$, $\theta_{e_U}(v)$ is an element of
the fiber of $p_{e_U}(x)$, where $p_{e_U}:e_U\rightarrow U$ is the
canonical projection.

 Let $A$ be an element of ${\cal L}'$, and $\bar A$ the projection
 of $A$ in ${\cal L}$ by the canonical map $\bar p:{\cal L}'\rightarrow
 {\cal L}'/{\cal H}={\cal L}$. Denote by $A^*$ the fundamental
 vector field generated by $A$ on $e_U$. We assume that ${\bar
 p}(\theta_{e_U}(A^*))={\bar A}$.

 Let ${\cal H}_{e_U}^*$ be the vector space of fundamental vectors
 generated by elements of ${\cal H}$. We assume that
 $\theta_{e_U}$ preserves ${\cal H}^*_{e_U}$, and its restriction
 to it is a projection.

 Let $h:e_U\rightarrow e_{U'}$ be a morphism in $P$, we assume
 that $h^*(\theta_{e_{U'}})=\theta_{e_U}$.

 \medskip

An horizontal path  in $e_U$ is a differentiable  path
$c:I\rightarrow e_U$ such that for each $t$ in $I$, the tangent
vector to the curve $c'(t)$ at $t$ is an element of ${\cal
D}_{e_U}$.

 The holonomy of a transitive distribution can be
defined as is defined the holonomy of a connection (see
Lichnerowicz [19] p. 62, Molino [29] p. 181): Let $x$ be an
element of $e_U$, the holonomy group $H_x$ at $x$, is the set of
elements $l'$ in $L'$ such that there exists an horizontal path
between $x$ and $x{l'}^{-1}$, of course if we replace $x$ by $hx$,
$H_{hx}=Ad(h^{-1})H_x$.

\medskip

Let $x$ be an element of $U$, and $c:I\rightarrow U$ a
differentiable path such that $c(0)=c(1)=x$. Consider $y$ an
element of the fiber of $x$. Since the distribution is transitive,
there exists an horizontal path over $c$ in $e_U$, $d:I\rightarrow
e_U$ such that $d(0)=y$. This is implied by the fact that a
transitive invariant distribution contains always a connection.
(See Molino [29] p. 181). The element $y{d(1)}^{-1}$ does not
depends of $y$ (compare with Lichnerowicz p. 94). It is called the
holonomy around $c$. The holonomy group $H^{e_U}_x$ at $x$ is the
set whose elements are holonomy around loops at $x$. The holonomy
group depends of the object since two objects of $P_U$ are not
always isomorphic.

Suppose that $U$ is contractible, then the holonomy group does not
depends of the object since it is invariant by the gauge
transformations which preserve the connection since the extension
is central, and all the objects of $P_U$ are isomorphic. This last
group can be computed by the Ambrose-Singer theorem. (See Molino
[29] p. [183]).

\section{\bf Sequences of fibered categories  in differentiable
categories.}

One of the main motivation of the introduction of gerbes theory in
differential geometry is the geometric interpretation of
characteristic classes. Let $N$ be a manifold. A well-known result
identifies the $2$-dimensional  integral cohomology space
$H^2(M,{\Z})$ of $N$, with the space of isomorphic classes of
$U(1)$-bundles defined on $N$. This identification is one of the
main tool used in quantization in physics. String theory has
created the need of finding such an interpretation for higher
cohomology classes. The space of $3$-dimensional integral
cohomology classes is the classifying space of $U(1)$-gerbes (See
Brylinski [8] p. 200). The second Pontryagin class which is an
element of $H^4(N,{\Z})$ has been interpreted with $2$-gerbes by
Brylinski and McLaughlin (see [9] p. 625).  In this part, we are
going define and apply the theory and sequences of fibered
categories analog to the theory defined  in Tsemo [40]) to study
characteristic classes in differentiable categories.

\medskip

{\bf Definition 9.1.}

A $2$-sequence of fibered categories is defined by the following
data:

A fibered category  $ p:P\rightarrow C$ over the Grothendieck site
$C$, such that:

Let $U$ be an object of $C$, and $e_U$ an object of $P_U$. Recall
that $e_U$ is a differentiable manifold. There exists a
correspondence which assigns to $e_U$ a 2-category $Q_{e_U}$ (see
Benabou for the definition of a 2-category) whose objects are
gerbes defined on $e_U$.

Let $c:U\rightarrow U'$ a morphism of $C$, the restriction functor
is the pull-back of gerbes.

\medskip

There exists a covering $(U_i)_{i\in I}$, such that for every
objets $e_i$ and $e'_i$ of $P_{U_i}$, there exists an isomorphism
between the $2$-categories  $Q_{e_i}$ and $Q_{e'_i}$.

\medskip

The set of automorphisms of an object of $Q_{e_U}$ can be
identified with sections of a sheaf $L$ defined on $C$, and this
identification is natural in respect to morphisms between objects
and restrictions.

\medskip

Let $P"$ be the category whose objects are objects of $Q_{e_U}$,
$e_U$ in $P_U$, and $P'$ the category whose objects are open
subsets of the manifolds $e_U$.  If $e$ in $Q_{e_U}$ and $e'$ in
$Q_{e_{U'}}$ are objects of $P"$, there exists an open subset $V$
of $e_U$, (resp. $V'$ of ${e_{U'}}$)  a Lie group $H$ such that
$e\rightarrow V$ is a $H$-bundle (resp. $e'\rightarrow V'$ is a
$H$-bundle). A morphism $h':e\rightarrow e'$ in $P"$ is  a
morphism of $H$-bundles such that there exists $h:e_U\rightarrow
e_{U'}$ in $P$ such that $h(V)\subset V'$, and the following
square is commutative:

$$
\matrix{e & {\buildrel{h'}\over{\longrightarrow}} & e'\cr
\downarrow &&\downarrow\cr V &{\buildrel h\over{\longrightarrow}}
& V'}
$$

\medskip

Our descent condition is expressed by the fact that we assume that
the correspondence $P"\rightarrow P'$ which assigns to $e$,  the
open subset $V$ of $e_U$  is a fibered category.

\subsection{ Classification $4$-cocycles  and sequences of
$2$-fibered categories.}

Before to attach to a $2$-sequence of fibered categories a
cocycle, we describe an automorphism $h$ above the identity of a
$H$-gerbe $p:P\rightarrow N$ over a manifold $N$.

The automorphism $h$ is defined by a family of functors $h_U$ of
$P_U$ above the identity, where $U$ is an open subset in $N$, such
that if $V$ is a subset of $U$  the following square commutes:

$$
\matrix{P_U & {\buildrel{h_U}\over{\longrightarrow}} & P_U\cr
\downarrow r_{V,U} &&\downarrow r_{V,U}\cr P_V
&{\buildrel{h_V}\over{\longrightarrow}} & P_V}
$$

where $r_{V,U}:P_U\rightarrow P_V$ is the restriction map.

Let $(U_i)_{i\in I}$ be a good covering of $N$, we assume that an
object $e_i$ of $P_{U_i}$ is a trivial $H$-bundle. Since the
square:

$$
\matrix{P_{U_i} & {\buildrel{h_{U_i}}\over{\longrightarrow}} &
P_{U_i}\cr \downarrow r_{U_i\cap U_j,U_i} &&\downarrow r_{U_i\cap
U_j,U_i}\cr P_{U_i\cap U_j} &{\buildrel{h_{U_i\cap
U_j}}\over{\longrightarrow}} & P_{U_i\cap U_j}}
$$

The automorphism $h$ is described by a family of morphisms
$u_{ij}:U_i\cap U_j\times H\rightarrow U_i\cap U_j\times H$ such
that $u_{ij}^lu_{jl}^i=u_{il}^j$. Thus by a $H$-bundle.

\medskip

Now can describe the classifying $4$-cocycle:

We assume that $L$ is commutative. Let $(U_i)_{i\in I}$ be a good
cover of the site $(C,J)$, and $e_i$ and object of ${P}_{U_i}$, we
choose a gerbe $d_i$ in $Q_{e_i}$. Since ${P}_{U_i\times_CU_j}$ is
connected, there exists a morphism: $ u_{ij}:e^i_j\rightarrow
e^j_i $, and a map $u_{ij}^*:d_j^i\rightarrow d^j_i$. The
automorphism
$c^*_{ijl}=u^*_{ij}u_{jl}^*u^*_{li}:d_i^{jl}\rightarrow d_i^{jl}$
is not above the identity. But

$$
c_{ijlm}={u^*}^{ij}_{ml}{c}^*_{ijl}{u^*}^{ij}_{lm}\circ
{{c}^*_{ijm}}^{-1}\circ {c}^*_{ilm}\circ {{c}^*_{jlm}}^{-1}
$$

is a morphism above the identity that we identifies with a
$1$-form defined on $U_{ijlm}$ which takes its values in
$L(U_{ijlm})$. The Cech-DeRham isomorphism identifies this with a
$4$-cocycle which takes its values in $L$ if $C$ is a manifold,
since $L$ is assumed to be commutative.

\bigskip

Before to give examples, we are going to describe a weak version
of a $2$-sequence of fibered categories.

\medskip

{\bf Definition 9.1.2.}

A $2$-sequence of torsor/fibered categories, is a $2$-sequence of
fibered categories  where the sheaf of categories $P\rightarrow C$
is in fact a torsor.

\medskip

We can associate to a $2$-sequence of torsor/fibered categories a
$3$-cocycle as follows:

\medskip

Let $(U_i)_{i\in I}$ be a good cover of $C$, and $e_i$ the object
of $P_{U_i}$, and $d_i$ an object of $Q_{e_i}$, and
$u_{ij}:e^i_j\rightarrow e^j_i$. There exists a morphism
$u_{ij}^*:d_j^i\rightarrow d_i^j$.  The morphism

$$
c^*_{ijl}={u^*}^j_{li}{u^*}^l_{ij}{u^*}^l_{jl}
$$

is  above the identity. It is a $1$-form defined on $U_{ijl}$
which is $L(U_{ijl})$ valued. The Cech Derham isomorphism
identifies this cocycle with a $3$-form defined on $C$ which is
$L$-valued.

\bigskip

{\bf Examples.}

\medskip

Let $H$ be a compact simple Lie group. Consider a $H$-principal
bundle $p:P\rightarrow C$ over a differentiable category $C$. We
can define the following $2$-sequence of torsor/fibered
categories:

Let $U$ be an object of $C$, the object of $Q_{P_U}$ are
$U(1)$-gerbes which induces on each fiber of $e_U\rightarrow U$ a
gerbe isomorphic to the canonical $U(1)$-gerbe on $H$. This
example is defined when $N$ is a manifold by Brylinski and
McLaughlin [9] p. 625). In this situation the classifying cocycle
is an element of $H^3(N,U(1))=H^4(N,{\Z})$.

\medskip

This construction can also be applied to a  subcategory $C'_N$ of
differentiable category $C_N$ associated to a generalized orbifold
(see definition 3.2.1) where  there exists a simple and compact
Lie group $H$ such that every object of $C'_N$ is of the form
$(P,H,\phi_P)$. We can also define a $4$-integral class on the
space whose elements are closure of leaves of a foliation endowed
with a bundle like metric.

\medskip

Let $p:P\rightarrow N$ be a $H$-principal gerbe over the manifold
$N$. Without restricting the generality, we suppose that for every
open subset $U$ of $N$, and for every element $e_U$ in $P_U$,
$e_U$ is a principal $H$-bundle. We can construct the following
$2$-sequence of fibered categories:

 The objects of $Q_{e_U}$ are $U(1)$-gerbes which induces on
the fiber of $e_U\rightarrow U$ the canonical $U(1)$-gerbe defined
on $H$.

The classifying cocycle of this $2$-sequence of fibered categories
is a $4$-class in $H^4(N,U(1))$ which can be identified with a
$5$-class in $H^5(N,{\Z})$.

\bigskip

{\bf References.}

\bigskip

[1] C. Albert; Introduction \`a l'\'etude des vari\'et\'es lisses.
Th\`ese Universit\'e de Montpellier 1974.

\medskip

[2] C. Albert, P. Molino; Pseudo-groupes de Lie transitifs.
D\'epartement de Math\'ematiques Universit\'e de Montpellier.
1983.

\smallskip

[3] M. Artin, A. Grothendieck, J-L. Verdier; Th\'eorie des topos
et cohomologie \'etale des sch\'emas. S\'eminaire de G\'eom\'etrie
alg\'ebrique du Bois-Marie. 4, Vol 1.

\smallskip

[4] M. Audin; Op\'erations Hamiltoniennes de tores sur les
vari\'et\'es symplectiques. (Quelques m\'ethodes topologiques).
I.R.M.A 1989 Strasbourg.

\smallskip

[5] J. Benabou; Introduction to bicategories; Lectures Notes in
Mathematics. {\bf 40} 1-77.

\smallskip

[6] X. Buff, J. Fehrenbach, P. Lochak, L. Schneps, P. Vogel;
Espaces de modules des courbes, groupes modulaires et th\'eorie
des champs. Panoramas et Synth\`eses. Soci\'et\'e Mathematiques de
France. 1999. Paris.

\smallskip

[7] L. Breen;  and W. Messing; Differential geometry of gerbes.
Preprint available at arxiv.org, published at: Advances in
Mathematics. {\bf 1998} (2005) 732-846.

\smallskip

[8] J-L. Brylinski; Loop spaces characteristic classes and
geometric quantization. Progress in Math. 1993. New York.

\smallskip

[9] J-L Brylinski; D-A McLaughlin; The geometry of degree four
characteristic classes and of line bundles on loop spaces I. Duke
Math. Journal. {\bf 75} (1994) 603-638.

\smallskip

[10] J-L Brylinski; D-A McLaughlin; The geometry of degree four
characteristic classes and of line bundles on loop spaces II. Duke
Math. Journal. {\bf 83} (1996) 105-139.

\smallskip

[11] W. Chen, Y. Ruan; A new cohomology theory for orbifold.
available at arxiv.org, and published by: Commun. Math. Phys. {\bf
248} (2004) 1-31.

\smallskip

[12] C. Debord; Holonomy groupoids of singular foliations. J.
Differential Geom. {\bf 58} (2001)  467-500.

\smallskip

[13] J. Giraud; M\'ethode de la descente. M\'emoires de la
Soci\'et\'e Math\'ematiques de France. 1964 Paris.

\smallskip

[14] J. Giraud; Cohomologie non Ab\'elienne. Springer 1971
New-York.

\smallskip

[15] W. Goldman; Geometric structures and varieties of
representations. Contemporary Mathematics. {\bf 74} (1988)
169-198.

\smallskip

[16] S. Kobayashi K. Nomizu; Foundations of differential geometry.
John Wiley 1963 New-York.

\smallskip

[17] J-L. Koszul; Homologie et cohomologie des alg\`ebres de Lie.
Bulletin de la Soci\'et\'e Math\'ematiques de France. {\bf 78}
(1950) 65-127.

\smallskip

[18] J-L. Koszul; Sur certains groupes de transformations de Lie.
Colloque de G\'eom\'etrie Diff\'erentielle Strasbourg, 1953.

\smallskip

[19] A. Lichnerowicz; Th\'eorie globale des connexions et des
groupes d'holonomie. Edizioni Cremonese 1956 Rome.

\smallskip

[20] D. Long, W. Reid; All flat manifolds are cusps of hyperbolic
orbifolds. Algebraic and Geometric Topology. {\bf 2} (2002)
285-296.

\smallskip

[21] L. Lupercio, B. Uribe;  Loop Groupoids, Gerbes, and Twisted
Sectors on Orbifolds. Preprint available at arxiv.org.

\smallskip

[22] M. Mackaay, R. Picken;  Holonomy and parallel transport for
Abelian gerbes. Preprint available at arxiv.org and published at:
Advances in Mathematics. {\bf 170} (2002) 287-339.

\smallskip

[23] D. McDuff, D. Salamon; Introduction to symplectic topology.
Oxford Sciences Publications, 1998.

\smallskip

[24] D. McDuff; Enlarging the Hamiltonian group. Journal of
Symplectic Topology {\bf 3} (2005) 481-530.

\smallskip

[25] A. Medina; Autour des connexions localement plates sur les
groupes de Lie. Th\`ese Universit\'e de Montpellier 1979.

\smallskip

[26] A. Medina; Groupes de Lie munis de pseudo-m\'etriques
bi-invariantes. S\'eminaire de G\'eom\'etrie Diff\'erentielle de
Montpellier. (1981-1982).

\smallskip

[27] A. Medina, P. Revoy; Alg\`ebre de Lie et produit scalaire
invariant. Annales Scientifiques de l'E.N.S. {\bf 18} (1985)
553-561.

\smallskip

[28] J. Milnor; Curvature of left-invariant metrics on Lie groups;
Advances in Mathematics. {\bf 21} (1976) 283-329.

\smallskip

[29] P. Molino; Champs d'\'el\'ements de contact sur un espace
fibr\'e principal diff\'erentiable. Annales de l'Institut Fourier.
{\bf 14} (1964) 163-219.

\smallskip

[30] P. Molino; Sur la $G$-g\'eom\'etrie des vari\'et\'es
diff\'erentiables. Annali della Scuola Normale Superiore di Pisa.
{\bf 23}  (1969) 461-480.

\smallskip

[31] P. Molino; Sur la g\'eom\'etrie transverse des feuilletages.
Annales de l'Institut Fourier. {\bf 25}  2 (1975) 279-284.

\smallskip

[32] P. Molino; Le probl\`eme d'\'equivalence pour les
pseudo-groupes de Lie m\'ethodes intrins\`eques. Bulletin de la
Soci\'et\'e Math\'ematiques de France. {\bf 108} (1980) 95-111.

\smallskip

[33] P. Molino, M. Pierrot; Th\'eor\`emes de slice et holonomie
des feuilletages Riemanniens singuliers. Annales de l'Institut
Fourier. {\bf 37} 4 (1987) 207-223.

\smallskip

[34] M.N. Nguiffo-Boyom; Pseudo-groupes de Lie extremaux et
quelques applications topologiques. Th\`ese Universit\'e de Paris
XI 1977.

\smallskip

[35] M.S. Raghunathan; Discrete subgroups of Lie groups. Springer
(1972) New-York.

\smallskip

[36] J-G. Ratcliffe, S-T. Tschantz; Fibered orbifolds and
crystallographic groups. Preprint available at www.arxiv.org
arXiv:0804.0427.

\smallskip

[37] I. Satake; The Gauss-Bonnet theorem for $V$-manifolds. J.
Math. Soc. Japan {\bf 9} (1957) 464-492.

\smallskip

[38] J-P. Serre; Cohomologies des groupes discrets. Annals of
Mathematics Studies. {\bf 70} (1971) 77-169.

\smallskip

[39] A. Tsemo; D\'ecomposition des vari\'et\'es affines. Bulletin
des Sciences Mathematiques. {\bf 125} (2001) 71-83.

\smallskip

[40] A. Tsemo; Non Abelian cohomology: the point of view of gerbed
tower. African Diaspora Journal of Mathematics. {\bf 4} (2004)
67-85.

\smallskip

[41] A. Tsemo; Gerbes, $2$-gerbes and symplectic fibrations. Rocky
Mountain Journal of Mathematics. {\bf 38} (2008) 727-778.

\smallskip

[42] A. Tsemo; $G$-structures and gerbes. In preparation.

\end{document}